\documentclass[10pt]{article} 
\usepackage[utf8]{inputenc} 
\usepackage[polutonikogreek,english]{babel}
\usepackage[or]{teubner}
\usepackage{xpatch}
\xapptocmd{\greektext}{\edef~{\string~}}{}{}
\usepackage{amsfonts}
\usepackage{amsmath}
\usepackage{amsthm}
\usepackage{titlesec}

\titleformat{\section}  
  {\fontsize{13}{15}\bfseries} 
  {\thesection.} 
  {1em} 
  {} 
  [] 

\titleformat{\subsection}  
  {\fontsize{12}{14}\bfseries} 
  {\thesubsection.} 
  {1em} 
  {} 
  [] 

\titleformat{\subsubsection}  
  {\fontsize{11}{13}\bfseries} 
  {\thesubsubsection.} 
  {1em} 
  {} 
  [] 

\usepackage{hyperref}
\usepackage{geometry} 
\usepackage{graphicx}
\graphicspath{ {c:/users/fpc-dprotopapas/documents/latex/} }

\usepackage{booktabs} 
\usepackage{array} 
\usepackage{paralist} 
\usepackage{verbatim} 
\usepackage{subfig}
\usepackage{multirow}

\usepackage{fancyhdr} 
\usepackage[nottoc,notlof,notlot]{tocbibind} 
\usepackage[titles,subfigure]{tocloft}

\title{The mathematics 
of periodic anthyphairesis as a basis for the full understanding of Plato's philosophy
}
 \author{Stelios Negrepontis and Athanase Papadopoulos}

\date{\today}

\begin{document}

\newtheorem{proposition}{Proposition}[section]
\numberwithin{proposition}{subsection}
\newtheorem{lemma}[proposition]{Lemma}
\newtheorem*{proposition*}{Proposition}
\newtheorem{definition}[proposition]{Definition}
\newtheorem*{definition*}{Definition}

  \maketitle

 \noindent      {\bf Abstract.} Even though Plato’s philosophy in ancient times was always closely associated with mathematics, modern Platonic scholarship, during the last five centuries, has  moved steadily toward de-mathematization. The present work aims to outline a radical re-interpretation of Plato’s philosophy, according to which the Platonic Idea, that is, the intelligible Being, has the structure of the philosophical analogue of a
  geometric dyad in a philosophic anthyphaeresis --- the precursor of modern continued
  fractions --- which was studied by the Pythagoreans, Theodorus and Theaetetus in
  relation with the discoveries of quadratic incommensurabilities. This mathematical
  structure is clearly visible in the Platonic method of Division and Collection, equivalently
  Name and Logos, equivalently True Opinion plus Logos, in the dialogues Theaetetus,
 Sophist, Statesman, Meno, and Parmenides. Equipped with this structure of an
 intelligible Being, we provide definitive answers to fundamental questions, that
 were not be resolved by Platonists, concerning the following topics: the dialectic
numbers, which are based on the anthyphairetic periodicity and the plus one rule,
stating that the dialectic number of terms of a sequence is the (number of) ratios  
   of successive terms plus one (stated in the Parmenides 148d-149d); the description of
 the intelligible being as an Indivisible Line, a statement bordering on the contradictory;
 the also seemingly contradictory Sophist ’s statement that “the not-Being is a Being”,
 based on the equalization of the two elements of the dyad defining an intelligible
 Being; the more general self-similar Oneness of an intelligible Being, based on the
equalization of all parts generated by the anthyphairetic division of an intelligible
Being; and finally the Third Man Argument in the Introduction to the Parmenides,
appearing as a threat for Plato’s theory, but essentially innocuous because of the
 self-similar Oneness. The third part of our study aims to prove that, contrary to
 the presently dominant interpretation of Zeno’s arguments and paradoxes as being
devoid of mathematical content, the analysis of Zeno’s presence in the Parmenides,
Sophist (via the Eleatic Stranger), and Zeno’s verbatim Fragments preserved by
Simplicius, show that Plato’s intelligible Beings essentially coincide with Zeno’s true
Beings, and hence that Zeno’s philosophical thought was already anthyphairetic, and
hence heavily influenced by the Pythagorean’s Mathematics. These findings run
 against Burkert’s claim that “ontology is prior to mathematics”. Modern Platonists
  have never obtained a clear description of the structure of an intelligible Idea in terms
  of the mathematics of periodic anthyphairesis, and thus were not able to answer
  fundamental questions, nor to realize the close connection of Zeno’s intelligible beings
  with Zeno’s true Beings.
 
\bigskip


\bigskip


\bigskip

\vfill\eject
  
  \section{Introduction}

 Plato's philosophy posits the existence of intelligible Beings, the Platonic Ideas, distinct and separate from the sensible entities.  An intelligible Being is an unchanging One, a monad, while in contrast, a sensible entity is always changing, consisting of many parts, coming to be and being destroyed. The intelligible Beings are superior to the sensible entities, but the sensibles participate in the world of Ideas.  A sensible entity has property A only by participating in the Idea of A. For instance, something is beautiful if it participates in the Idea of beauty. In Plato's philosophy, knowledge is knowledge of Ideas, for example, of the Idea of beauty, in contrast to opinion, which is about sensibles; for example, opinion may be about beautiful things. The difference between knowledge and opinion is essentially the difference between the world of Ideas and the world of the sensibles. To understand this distinction by Plato, mathematics is always in the background. Let us quote a passage from the \emph{Philebus} (51b9-c7): 
\begin{quote}\small
For when I say beauty of form, I am trying to express, not what most people would understand by the words, such as the beauty of animals or of paintings, but I mean, says the argument, the straight line and the circle and the plane and solid figures formed from these by turning-lathes and rulers and patterns of angles; perhaps you understand. For I assert that the beauty of these is not relative, like that of other things, but they are always absolutely beautiful by nature.
\end{quote}

Plato's writings have been studied incessantly for the last five centuries, that is, after they were translated from the Greek into Latin in the ambitious program undertaken by Marsilio Ficino, in 1484, under the auspices of the Florentine Medici.
The influence of these writings on  modern philosophy and more generally on human thought cannot be overestimated.
We adhere to what  Alfred North Whitehead wrote in 1929 \cite[p. 39]{Whitehead}:  

\begin{quote}\small
 The safest general characterization of the European philosophical tradition is that 
it consists of a series of footnotes to Plato. 
I do not mean the systematic scheme of thought 
which scholars have doubtfully extracted from his writings. 
I allude to the wealth of general ideas scattered through them.
\end{quote}

It is also worth remembering that at the school he founded, the \emph{Academy}, Plato  was essentially teaching Mathematics. We learn this for example from a passage by Aristoxenus of Tarentum, the famous peripatetic fourth century BC music theorist, author of the oldest almost complete treatise on music, the \emph{Elementa Harmonica}. In this treatise, Aristoxenus parenthetically reports on Plato's teaching.
We read there:

\begin{quote}\small

  Such was the condition, as Aristotle used often to relate, of most of the audience that attended Plato's lectures on the Good. 
They came, he used to say, every one of them, in the conviction that they would get from the lectures some one or other of the things that the world calls good ; riches or health, or strength, in fine, some extraordinary gift of fortune. 
But when they found that Plato's reasonings [logoi] were 
of sciences [mathemata/mathematics] and numbers, and geometry, and astronomy, and of good and unity as predicates of the finite, 
[\ldots]  their disenchantment [paradoxon] was complete. 
The result was that some of them sneered at the thing, while others vilified it (Aristoxenus, \emph{Elementa Harmonica} 39,8-40,4).
\end{quote}

 Plato's high regard for Geometry and its relation to knowledge
 can be seen in the following passage from the \emph{Republic} (527b3-c3):
 
  \begin{quote}\small

Socrates: And must we not agree on a further point?

Glaucon: What?

 Socrates:  That it is the knowledge of that which always is, 
and not of a something which at some time comes into being and passes away.

Glaucon:  That is readily admitted. For geometry is the knowledge of the eternally existent.

Socrates: Then, my good friend, it would tend to draw the soul to truth, 
and would be productive of a philosophic attitude of mind,    
directing upward the faculties that now wrongly are turned earthward.

Glaucon:   Nothing is surer.
  
 Socrates: Then nothing is surer than that 
we must require that the men of your Fair City shall never neglect geometry, for even the by-products of such study are not slight."

\end{quote}

The fact that Plato was a mathematician is confirmed by Proclus of Lycia, the Greek Neoplatonist philosopher and one of the major philosophers of late antiquity, whose opinion is generally considered to be faithful.
In his  \emph{Commentary to Euclid}, 
Proclus includes a Synopsis of Greek Mathematics up to Euclid, information whose origin is generally considered to be due Eudemus, Aristotle's student. In this Synopsis, Proclus includes Plato among the mathematicians. This is what he writes about him:

\begin{quote}\small
[\ldots] Plato, who appeared after them, greatly advanced mathematics in general and
geometry in particular because of his zeal for these studies. It is well known that
his writings are thickly sprinkled with mathematical ratios (mathematikois logois) and that he everywhere tries to arouse admiration for mathematics among students of philosophy.\footnote{This translation is by Glen Morrow, who however obscures the meaning of Proclus' passage by rendering the Greek words ``mathematikois logois", possessing the clear meaning ``mathematical ratios", by the expression ``mathematical terms", neutralizing its real meaning. We have modified the translation accordingly.}

\end{quote}

Even though Plato's philosophy in ancient times was always closely associated with mathematics, modern Platonic scholarship, during the last five centuries, has moved steadily toward de-mathematization. 
 Bertrand Russell, one of the main founders of contemporary logic and a central figure of Analytic philosophy, wrote in 1945 \cite[p. 132]{Russell}: 

\begin{quote}\small
It is noteworthy that 
modern Platonists, almost without exception, are ignorant of mathematics,
in spite of 
the immense importance that Plato attached to arithmetic and geometry, and 
the immense influence that they had on his philosophy. 
This is an example of the evils of specialization: 
a man must not write on Plato unless he has spent so much of his youth on Greek
as to have had no time for the things that Plato thought important.
\end{quote}

Even though Russell did not pinpoint the exact nature of the mathematics needed to understand Plato's philosophy, nevertheless he located the substantial difficulty. One might think that some Platonists would take Russell's diction seriously and try, even with the delay of centuries, to discover the as yet unknown mathematical base of Plato's philosophy. But in reality, during the interval from 1946 to our days, Platonists have distanced themselves even further from Mathematics.

 Although the whole work of Plato plus the numerous scholia on it have been preserved, and even if this work has been studied most intensively by able and talented scholars for the last five centuries, nevertheless there is a long list of fundamental questions that have not been answered in any satisfactory way but instead remain still open. 
 
According to Plato,  Geometry leads to the knowledge of the Intelligible Beings.  But the high regard of Plato for Geometry does not apply to geometers in general, precisely because, according to him, they relied on axioms, postulates, diagrams etc., instead of using the method of  anthyphairesis  (the precursor of modern continued fractions),  which was studied by the Pythagoreans, Theodorus and Theaetetus in relation with the discoveries of quadratic incommensurabilities and is the main object of interest in the following sections. This has been commented on in \cite{N-geometry}. To the best of our knowledge, the following are the scholars that have detected
some presence of anthyphairesis in Plato's philosophy: (i) A. E. Taylor
\cite{Taylor}, (ii) D'Arcy W. Thompson \cite{Thompson}, (iii) C. Mugler \cite{Mugler}, (iv) J. Vuillemin
\cite{Vuillemin}, (v) D. Fowler \cite{Fowler}.

 The present work aims to outline a radical re-interpretation of Plato's philosophy, according to which the Platonic Idea, that is, the intelligible Being, has the structure of the philosophical analogue of a geometric dyad in a philosophic anthyphairesis.
 By analyzing and revealing the underlying mathematical structure of his philosophy, we are able to provide definite answers to some fundamental questions that are at the basis of this philosophy. Here is a list of some of these questions, followed by references to works done on them by the first
author, which started in 1996 with \cite{N1997, N1999}, together with collaborators.
\begin{enumerate}

\item What are the Mathematics that are crucial and fundamental for Plato's philosophy?
\cite{N2012}

\item  What is the role of incommensurability in Plato philosophy?
\cite{NFB}

\item How can we interpretat  the statements, in \emph{Theaetetus, Sophist, Meno} and \emph{Parmenides}, that the knowledge of an intelligible Being is True Opinion plus Logos, or Name plus Logos? What is the Platonic meaning of Logos, and its importance in Plato's philosophy?
 \cite{N2012, N2018, N2026}  

\item  What is the structure and the nature of Plato's intelligible Beings and of the Platonic Ideas?
What is the Platonic meaning of Logos?
  \cite{N2012, N2018, N2024a, N2026}  
 \item What is the relation between the Platonic Ideas and the sensibles? 
what is the meaning of participation of the sensible in an intelligible?
 \cite{N-geometry, N2024a}
 \item What is the definition and role of the dialectic/Platonic numbers?
 \cite{N2026}  
 \item What is the nature of the Indivisible Line?
  \cite{N-geometry, N2024b}
 \item In what sense is the Being equalized with the One in the second hypothesis of the \emph{Parmenides}, and in what sense is every intelligible Being a One?
 \cite{N2026} 
 \item How can we explain Plato's statement in the \emph{Sophist} that the not-Being is a Being?
 \cite{N2026} 
 \item What is the meaning of the definition of an intelligible Being in the \emph{Sophist} 247ff?
 \cite{N2026} 
 \item How is the Third Man Argument in the Introduction of Plato's \emph{Parmenides} resolved?
  \cite{N2026} 
 \item In what sense does an intelligible Being satisfy the compresence of almost contradictory properties, such as Infinite and Finite, One and Many, in Motion and at Rest?
  \cite{N-Zeno-2020, N-Zeno} 
 \item How is the mystery of the Receptacle in the \emph{Timaeus} explained?
  \cite{N-2005, NP2024}
 \item What is the role of Zeno in Plato's philosophy?
What is the relation between Plato's intelligible Beings and Zeno's true Beings?
  \cite{N-Zeno-2020, N-Zeno} 
 \item What is the purpose and interpretation of Zeno's arguments and paradoxes?
  \cite{N-Zeno-2020, N-Zeno} 
 \item  What is the relation of Zeno arguments and paradoxes with the Pythagoreans?
\cite{N-Zeno-2020, N-Zeno}
 \item What is the status of Burkert's thesis that ``philosophy is prior to mathematics" and of Netz's assessment that ``Pythagoras the mathematician finally perished in 1962 AD"?
  \cite{N-Zeno-2020, N-Zeno, NF2025}
 \item How to explain the tremendous influence that Plato's philosophy had on Western thought?
\end{enumerate}

 Ending this introduction by a quote from a modern mathematician on Plato, we reproduce a sentence from  Bourbaki's \emph{\'Eléments d'histoire des mathématiques}\cite[p. 12]{Bourbaki-H}: 
\begin{quote}\small

It has been said that Plato was almost obsessed with mathematics;
without being himself an inventor in this field himself, he became, from a certain point in his life, acquainted with the discoveries of contemporary mathematicians (many of whom were his friends or students), and never ceased to take a direct interest in them, even going so far as to suggest new directions for research;
thus, in his writing, mathematics constantly serves as an illustration or a model (and sometimes even nurtures, as with the Pythagoreans, his 
inclination for mysticism).

\end{quote}

\section{The Mathematics of periodic anthyphairesis}\label{s:Periodic-anthyphairesis}
 
 In this section, we review the mathematical notion of anthyphairesis and the results we need.

We start by recalling the \emph{anthyphairesis} of two natural numbers. This is essentially the Euclidean algorithm, introduced in
 Propositions VII.1,2  of the \emph{Elements}.

Let $a>b$ be two natural numbers, and suppose that

$a= k_0b+c_1$, with $b>c_1$,

$b= k_1c_1+c_2$,  with $c_1>c_2$,

$c_1=k_2c_2+c_3$, with $c_2>c_3$,

etc.

By the so-called Principle of the Least (which is equivalent to the Principle of Mathematical Induction), this process of anthyphairesis of two natural numbers is always finite. Thus the two last steps in this process will necessarily have the following form: 

$c_{n-2}=k_{n-1} c_{n-1}+c_n$, with $c_{n-1}>c_n$,

  $c_{n-1}=k_nc_n$.

The number we get, $c_n$, is the Greatest Common Divisor of $a$ and $b$. 

Now we pass to \emph{geometric anthyphairesis}. This is dealt with in  Proposition X.2 of the \emph{Elements}.
  It concerns magnitudes (lines, areas, volumes, weight, time, etc.), that is, objects that are more general than natural numbers.
 
 Let $a, b$ be two magnitudes, with $a>b$.

 The \emph{anthyphairesis of $a$ to $b$} is the following finite or infinite sequence of  divisions:

$a   = k_0 b+ e_1$,           with $b>e_1$,

$b   = k_1e_1 +e_2$,           with $e_1>e_2$,

$\ldots$

$e_{n-1}= k_n e_n+ e_{n+1},$       with $e_n >e_{n+1}$,

$e_n    = k_{n+1} e_{n+1}+ e_{n+2}$, with $e_{n+1}>e_{n+2}$,

\ldots 

The \emph{quotients of the anthyphairesis of $a$ to $b$} is the sequence of successive natural numbers $k_0, k_1,\ldots, k_n, k_{n+1},\ldots$. We use the following notation:
 
\[\mathrm{Anth}(a, b)=[k_0, k_1,\ldots, k_n, k_{n+1},\ldots].\]

\medskip

  Proposition X.2 of the \emph{Elements} gives an
 \emph{Anthyphairetic Criterion for incommensurability} which we state as follows:
 \emph{ Two magnitudes $a, b$, with $a>b$ are incommensurable 
if (and only if) the anthyphairesis of $a$ to $b$ is infinite.}

\medskip

 Euclid's proof of X.2 is unnecessarily based on Eudoxus' Condition (Definition V.4 of the \emph{Elements}). An older proof, elementary, with no reliance on Eudoxus' Condition is already suggested in the passage  158b-159a of Aristotle's \emph{Topics} where the latter refers to a pre-Eudoxean theory of ratios of magnitudes, namely to a theory before the one presented in Book V of the \emph{Elements}. This theory, doubtlessly due to Theaetetus, is based on the following

\medskip

\noindent {\it Definition.} If $a, b, c, d$ are four magnitudes, with $a>b$, $c>d$, 
then $a/b=c/d$ if and only if $\mathrm{Anth}(a,b)=\mathrm{Anth}(c,d)$.
 A reconstruction of Theaetetus' theory of ratios of magnitudes has been given in \cite{NP2023, NP2025}.
  
\medskip

A trivial remark is that if a tail of the sequence of two magnitudes coincides with the sequence itself, then the sequence is periodic
An immediate consequence of the anthyphairetic definition of proportion of magnitudes and of this trivial remark is the following fundamental
 
 \bigskip
 
 \noindent {\it Proposition  (Logos Criterion for the periodicity of a geometric anthyphairesis)}.
 Let $a, b$ be two magnitudes, $a>b$, with anthyphairesis

$$a  = k_0 b+ e_1,                              \ \mathrm{with} \   b>e_1,$$  
$$b   = k_1e_1 +e_2,                              \ \mathrm{with} \   e_1>e_2,$$
$$\ldots$$
$$e_{n-1}    = k_n e_n+ e_{n+1},                      \ \mathrm{with} \  e_n >e_{n+1},$$
$$
e_n    = k_{n+1} e_{n+1}+ e_{n+2},               \     \mathrm{with} \  e_{n+1}>e_{n+2},
$$
$$ \ldots  $$ 
$$ e_{m-1}    = k_m e_m+ e_{m+1}, \                   \mathrm{with} \  e_m >e_{m+1},$$
$$
e_m    = k_{m+1} e_{m+1}+ e_{m+2},                 \     \mathrm{with} \  e_{m+1}>e_{m+2},
$$
     $$   \ldots, $$

such that for some indices $n<m$ we have  

        $$e_n/e_{n+1}= e_m/e_{m+1}             \ \  \mathrm{(Logos \  Criterion)}. $$

Then the anthyphairesis of $a$ to $b$ is eventually periodic. In fact,

$\mathrm{Anth}(a,b)= [k_0, k_1,\ldots, k_m, \ \mathrm
{period} \ (k_{m+1}, k_{m+2},…, k_n)]$.

\bigskip

\noindent \emph{Proof}.  

$\mathrm{Anth}(a,b)=[k_0, k_1,\ldots,  k_m, k_{m+1}, k_{m+2},\ldots, k_n,\ldots]$

$=[k_0, k_1, \ldots, k_m, k_{m+1}, k_{m+2},\ldots, k_n, \mathrm{Anth}(c_n, c_{n+1})]$

$=[k_0, k_1,\ldots, k_m, k_{m+1}, k_{m+2},\ldots, k_n, \mathrm{Anth}(c_m, c_{m+1})]$

$=[k_0, k_1,\ldots, k_m, k_{m+1}, k_{m+2},\ldots, k_m, k_{m+1},k_{m+2},\ldots, k_n,\mathrm{Anth}(c_n, c_{n+1})]$

\ldots

$=[k_0,k_1,\ldots, k_m, period (km+1,km+2,…, kn)].$

\begin{figure}[!h] 
\centering
\includegraphics[width=0.35\linewidth]{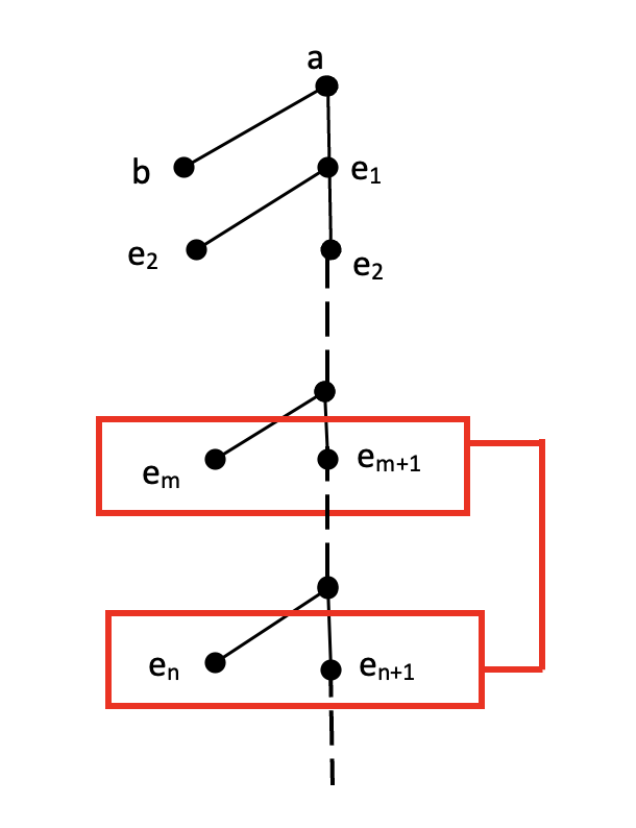}
\caption{\small {Abbreviated representation of the Logos Criterion}}
\label{fig:Anthyphairesis6}
\end{figure}

In the rest of this paper, we shall use several times the representation of the Logos Criterion given in Figure \ref{fig:Anthyphairesis6}.

 \bigskip

To conclude this section, let us quote again a passage from Bourbaki's \emph{Histoire des mathématiques} in which the author attributes the invention of anthyphairesis to the Pythagoreans \cite[p. 110]{Bourbaki-H}: 

\begin{quote}\small If $a_1$ and $a_2$ are two integers such that $a_1\geq a_2$, we define inductively $a_n$
(for $n\geq 3$) as the remainder of the Euclidean division of $a_{n-2}$ by $a_{n-1}$; 
if $m$ is the smallest index such that $a_m=0$,  then $a_{m-1}$ is the GCD of $a_1$ and $a_2$. This is
the transposition into the domain of integers of the method of successive subtractions
(sometimes also called anthyphairesis) for finding the common
measure of two quantities. This method undoubtedly dates back to the Pythagoreans and
seems to have been the basis of a pre-Eudoxian theory of irrational numbers.
\end{quote}

 Bourbaki did not have the arguments about the Pythagorean origin and dating of the discovery of the notion of anthyphairesis that we shall see all throughout this essay, 
nevertheless he had the right intuition.

\medskip

   In the next four sections we exhibit the structure of Plato's intelligible Being as a dyad in the philosophic analogue of periodic anthyphairesis, a structure clearly visible in the dialogues \emph{Theaetetus, Sophist, Meno} and (although much less so), in the \emph{Parmenides}, with a crucial role for Logos in this dyad, this Logos being the analogue of the Logos Criterion for periodicity. 
   
   \section{Plato's \emph{Theaetetus}} \label{Theaetetus}
  In the \emph{Theaetetus},  Plato prepares the ground for a double illumination, which is at the same time mathematical, on the nature of Theaetetus' mathematical achievements, and philosophical, on the geometrical nature of Plato's intelligible Beings.
Here are the four steps setting up this preparation:

(i) A \emph{Question} by Socrates at 145d4-e9: What is it that constitutes philosophical knowledge of the intelligible Beings?

(ii) An \emph{Answer} by Theaetetus at 147c7-148b2, connecting philosophy to mathema\-tics: The question on philosophical knowledge is similar to the one that Theaetetus discovered on geometrical knowledge when he proved the geometrical incommensurability of a power $a$ to line $b$, such that $a^2=Nb^2$, with $N$ being a non-square natural number.

(iii) The \emph{Urge} by Socrates for Imitation at 148c9-d7 (``peiro mimoumenos", 148d4-5): ``Try to imitate the geometrical answer for the question on philosophical knowledge."

(iv) The \emph{Crucial role of Logos for Knowledge} at 201d-202c: The knowledge of an intelligible Being is described as ``True Opinion plus Logos", equivalently as ``Name plus Logos".

Plato calls the method for acquiring the knowledge of an intelligible Being as Name and Logos, equivalently as True Opinion and Logos, equivalently as Division and Collection (cf. \emph{Phaedrus} 265d3-266c1). As we shall see, the three terms Name, True Opinion, and Division refer to an initial segment of the anthyphairetic division, whereas Logos and Collection refer to the Logos Criterion establishing anthyphairetic periodicity, by which all the infinite multitude of parts generated by the anthyphairetic division are collected into one, as indeed accomplished, by the anthyphairetic periodicity, cf. Section \ref{s:dialectic} below. In this way, the mathematical achievements  of Theaetetus are intimately connected with the mathematics of periodic anthyphairesis (\cite{NP2023, NFB2024}).

\section{The Knowledge of the intelligible being ``the Angler" in the \emph{Sophist}}\label{s:Angler}

In this section, we show that the Knowledge of the intelligible Being ``The Angler" as Name + Logos in the \emph{Sophist} 218b-221c has the structure of the philosophic analogue of periodic anthyphairesis, with Logos corresponding to the Logos Criterion for periodicity.
The method of ``Name and Logos" (cf. \emph{Theaetetus} 201e2-202b5, \emph{Sophist} 218c1-5, 221a7-b2, 268c5-d5) is exemplified at the beginning of the \emph{Sophist} 218b-221c by the definition of the Being Angler.

\begin{figure}[!h] 
\centering
\includegraphics[width=0.7\linewidth]{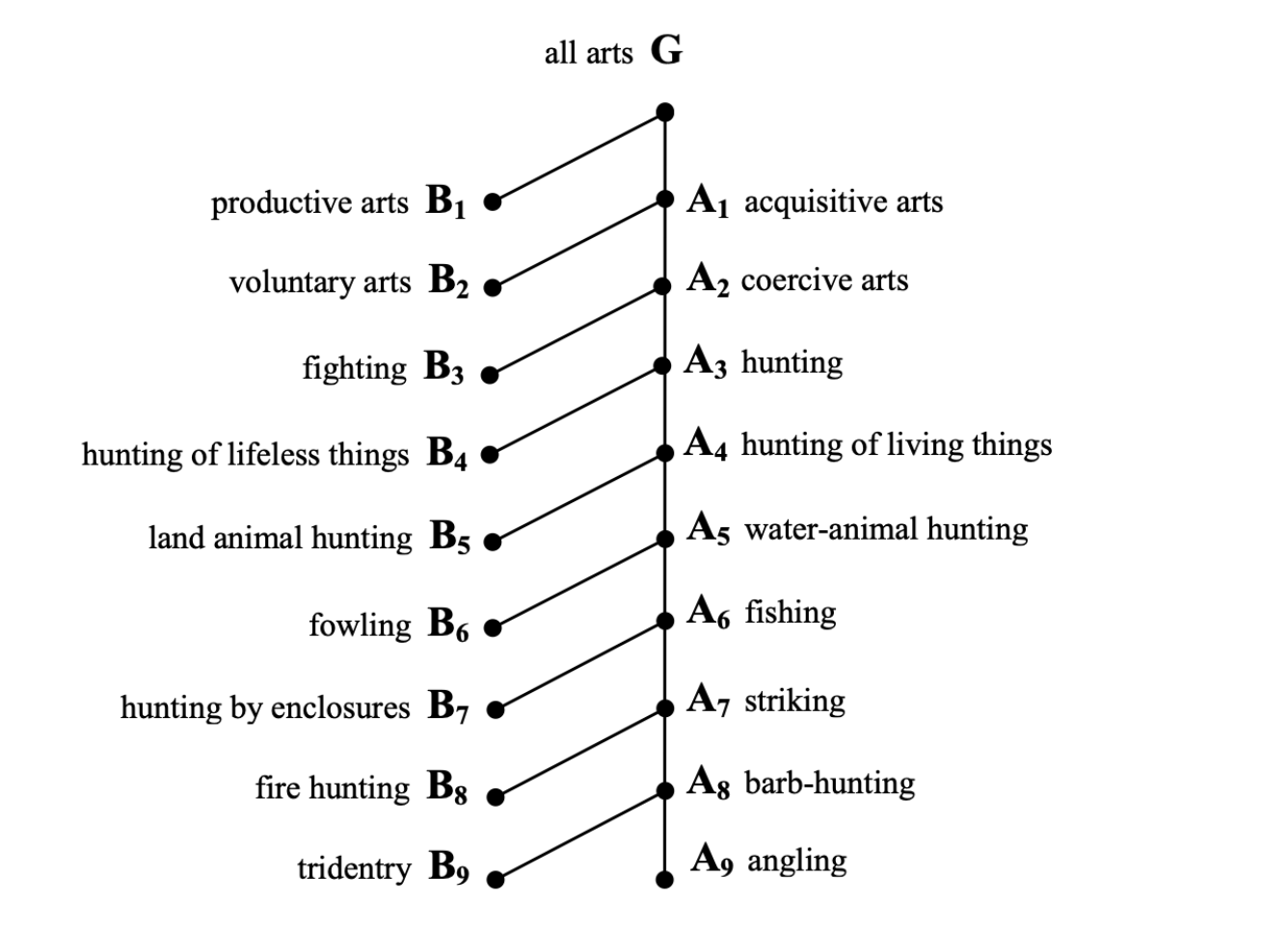}
\caption{\small {The Division of the Angler (\emph{Sophist} 218b-221c)}}
\label{fig:Sophist1}
\end{figure}

In Fig. \ref{fig:Sophist1}, we have reproduced the binary division process by which Plato leads to the definition of the Angler in the \emph{Sophist} 218b-221c.  The crucial remark by Plato, which gives the Logos Criterion, is the one in \emph{Sophist} 220e2-221b1, saying that the ratio of 
Tridentry B9 to Angling A9 
``is equal,
from above downward (‘anothen eis to kato'), to from below upwards (‘katothen eis tounantion ano')".

 In the next passage, 
 \emph{Sophist} 220e2-221b1, it is also stated that the ratio
B6/A6 is equal to
from above downward (``anothen eis to kato") to from below upwards (``katothen eis tounantion ano").
The complete passage from the \emph{Sophist} is the following:
\begin{quote}\small
but we have acquired also a satisfactory ‘Logos' of the thing itself.
For (γὰρ) 
of art as a whole, half was acquisitive,
and of the acquisitive, half was coercive, 
and of the coercive, half was hunting,
and of hunting, half was animal hunting,
and of animal hunting, half was water hunting, 
and of water hunting [A5],
the whole part from below (τὸ κάτωθεν τμῆμα ὅλον) was fishing [A6],
and of fishing, half was striking,
and of striking, half was barb-hunting [A8],
and of this [A8]
the part in which the blow is pulled 
from below upwards (τὸ περὶ τὴν κάτωθεν ἄνω) was angling [A9]. (221b1-c3). 

\end{quote}

\begin{figure}[!h] 
\centering
\includegraphics[width=1\linewidth]{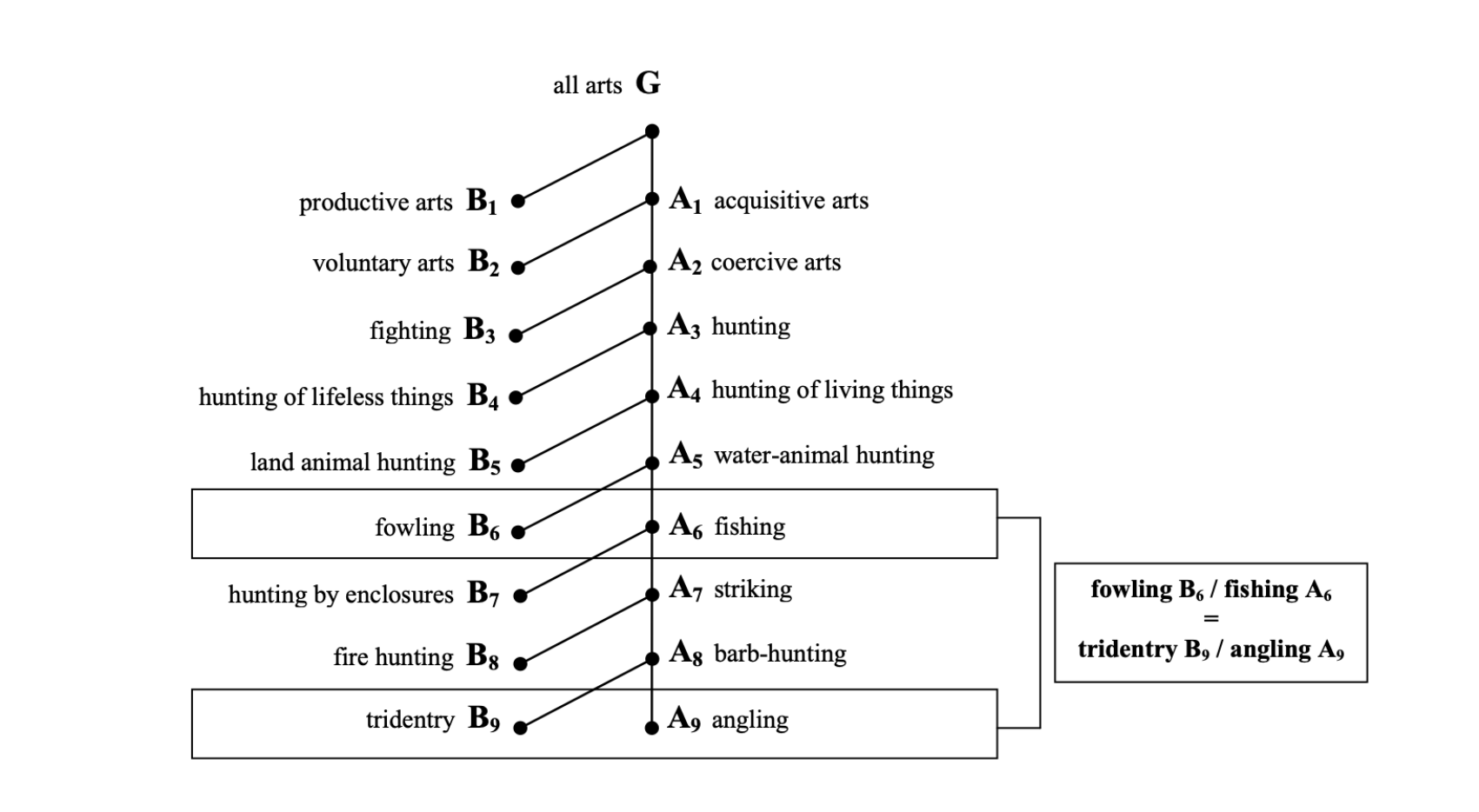}
\caption{\small {The equality expressing the Logos Criterion for anthyphairetic periodicity in the definition of the Angler}}
\label{fig:Sophist2}
\end{figure}

From this, we deduce  a satisfactory ``Logos" of the thing itself (221a7-b2),
namely the Logos:

\centerline{Fowling B6 / Fishing A6= Tridentry B9 / Angling A9,}

\noindent which we see a philosophical analogue of the geometrical Logos Criterion for anthyphairetic periodicity. We have represented this schematically in Figure \ref{fig:Sophist2}.

\section{The Knowledge of the diameter to the side of a square in the \emph{Meno}} 

In this section, we show how the knowledge of the diameter to the side of a square is described as True Opinion plus Logos in the \emph{Meno} 80d-86e, 97a-98b, where Logos corresponds again to the Logos Criterion of anthyphairetic periodicity.

 In the \emph{Meno}, Plato treats the diameter $a$ to the side $b$ of a square as an intelligible Being, and thus provides an ideal link between Geometry and Plato's Philosophy. 
The Knowledge of the diameter $a$ to the side $b$  is given by True Opinion plus Logos.
In this setting, True Opinion consists of the first two steps of the anthyphairesis of $a$ to $b$  (\emph{Meno} 82e2-83e10):
\begin{equation}\label{e:8.2}
 a=b+c_1  \ \ \mathrm{and}   \ \ 
 b=2c_1+c_2.
\end{equation}

Indeed, by the (isosceles) Pythagorean theorem, a line $a$ constructed as the diameter of the square with given side $b$ satisfies $a^2=2b^2$; then,
$b^2<a^2=2b^2<4b^2$, hence  $b<a<2b$, hence $a=b+c_1$, $b>c_1$ and
$9b^2>8b^2$, therefore $9b^2>4a^2$, or $3b>2a$, hence $b>2(a-b)$, hence $b>2c_1$, hence $b=2c^1+c^2$.

 As mentioned in \S \ref{s:Periodic-anthyphairesis}, the Theaetetean definition of proportion of magnitudes is 
$a/b=c/d$ if $\mathrm{Anth}(a,b)=\mathrm{Anth}(c,d)$.
In the present situation,

 \begin{equation} \label{e:8.4}
 \mathrm{Logos \  is \ the \ ratio} \ b/c_1=c_1/c2.
 \end{equation}
  Indeed, we have 
$c_1=a-b$, hence $c_2=b-2c_1=b-2(a-b)=3b-2a$, hence 
$b\cdot c_2=b(3b-2a)=3b_2-2ab$, and
$c_1^2=(a-b)(a-b)=a^2+b^2-2ab=3b^2-2ab$.
Hence $b\cdot c^2=c_1^2$ and $b/c_1=c_1/c_2$. 
From this we can say exactly what is the Knowledge of the diameter (with respect to side), which is the full infinite anthyphairesis of the diameter $a$ to the side $b$, is $[1, 2,2,2,\ldots]$. Indeed, we have, 
\[\mathrm{Anth} (a, b) = [1, 2, \mathrm{Anth} (c_1, c_2)] \ \mathrm{(by \ Eq. \ } \ref{e:8.2})
 =[1, 2, \mathrm{Anth} (b, c_1)]
\]
\[= [1, 2, 2, \mathrm{Anth} (c_1, c_2)] \ \mathrm{(by \ Eq. \ } \ref{e:8.4})
=[1, 2,2, \mathrm{Anth} (b, c_1)]=
[1, 2, 2, 2, \mathrm{Anth} (c_1, c_2)]\]
\[=…
= [1, 2, 2, 2,…].\]

The fact that Plato in the \emph{Meno} treats the diameter $a$ to the side $b$ of a square as a philosophical intelligible Being, and calls 
the first two steps of the anthyphairesis True Opinion, 
the ratio $b/c_1=c_1/c_2$ Logos and 
the resulting full anthyphairesis of $a$ to $b$ Knowledge (episteme),
strongly suggests that
Plato's intelligible Being in general possesses the structure of a philosophical analogue of a dyad in periodic anthyphairesis, with
True Opinion being an initial finite segment of this anthyphairesis (taken up to a point that we can apply Logos, Plato's term for the Logos Criterion, cf. \ref{s:Periodic-anthyphairesis}).
Logos is a condition on the remainders appearing already in the True Opinion that ensures anthyphairetic periodicity. Thus, we have again, 

\centerline{Knowledge=True Opinion+Logos}

\noindent as in the \emph{Theaetetus} general plan (\S \ref{Theaetetus}).
 This is precisely the content of Plato's \emph{Meno} 80d-86e, 97a-98b, as shown in \cite{N2024a}.
The term ``recollection" by which Knowledge is described in the \emph{Meno} refers to the repetition of the same ratio in ordre to satisfy the Logos criterion for periodicity.

The interpretation of \emph{True Opinion} of an intelligible Being, found in the \emph{Meno}, as a finite initial segment of the infinite anthyphairesis of the intelligible Being $F$, is a valuable element for understanding the important passage \emph{Republic} 475e4-480a13. According to it, the knowledge of a sensible entity participating in an intelligible Being is identical to a true opinion of the intelligible Being. The model for this relation of participation of a sensible in an intelligible Being is the anthyphairesis of a side and diameter numbers in relation with the infinite anthyphairesis of the diameter to the side of a square. 
In the \emph{Timaeus}, Plato, in order to avoid a sensible form of the Third Man Argument, which is the subject of  \S \ref{s:Third} below, is forced to modify drastically the meaning of participation of a sensible into an intelligible and he introduces the Receptacle, by means of the four canonical solids that the Demiurge uses in the \emph{Timaeus}, see \cite{NP2024}.

The mathematical basis of the participation of the sensibles in the intelligible Beings, either in the initial sense or in the modified sense of the \emph{Timaeus}, has not been understood by modern Platonists.
 
\section{The intelligible Being ``the Sophist" forms the philosophical analogue of a dyad in periodic anthyphairesis}
In this section, we report on the Knowledge of the  Sophist by Name and Logos as the philosophic analogue of periodic anthyphairesis, a fact which we deduce from the passage 264b-268d, 258b5-7 of the \emph{Sophist}.

The definition of the intelligible Being ``The Sophist" in the \emph{Sophist} relies on a fundamental proportion of the Divided Line, on which Plato reports in the \emph{Republic} 509d1-510b1: 

\begin{quote}\small

-- Would you be willing to say, said I, that 
the division in respect of reality and truth or the opposite is expressed by the proportion:

\emph{as is the opinable [B] (to doxaston) to the knowable [A] (to gnoston),
so is the likeness [B2] (to homoiothen) to that [B1] of which it is a likeness (to hoi homoiothei)?}"

-- I certainly would. 
\end{quote}

We conclude that the ratio of Plato's Divided Line is the following:

\bigskip

  \centerline{Knowledge[A] /True Opinion[B]} 
   
  \centerline{=}
 
  \centerline{real things in themselves[B1] /images of the real things [B2].}

\bigskip

We now prepare the ground for the discovery of the Logos of the intelligible Being ``Sophist". 
We shall see that the fundamental proportion of the Divided Line in the \emph{Republic} 509d1-510b1 is essentially the Logos Criterion for the definition of the Sophist.

\begin{figure}[!h] 
\centering
\includegraphics[width=0.6\linewidth]{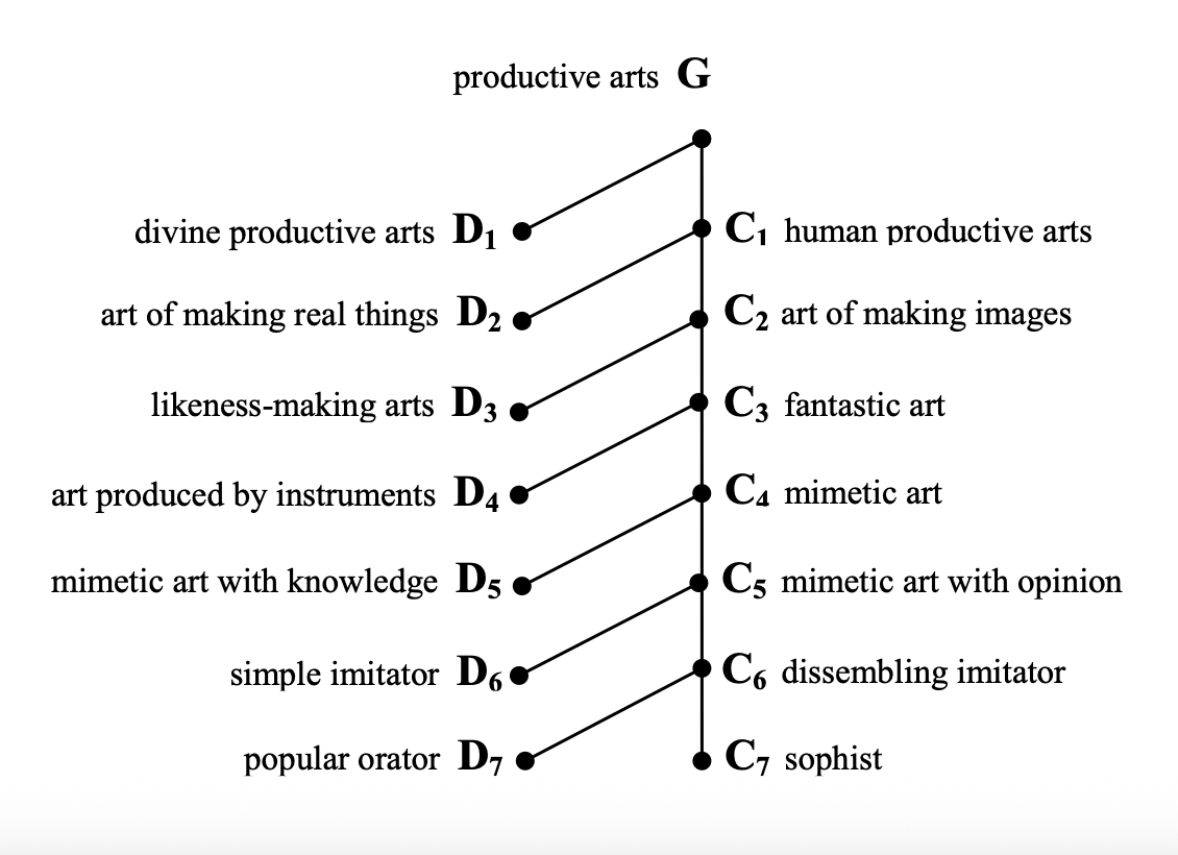}
\caption{\small {Division of the Sophist (\emph{Sophist} 264b-268d)}}
\label{fig:Sophist3}
\end{figure}

In the \emph{Sophist}, the Division/Name for the Sophist, starts with 
a Genus, in this case ``all the productive arts" and 
 proceeds by binary division of each Genus into two species, where the next Genus is that species of the previous step in the Division that contains the entity to be defined, in this case the Sophist, and ending with the division-step that produces the Sophist as a species. 
The whole division scheme is represented in Figure \ref{fig:Sophist3}.

  Using the ratio of the Divided Line, the Logos Criterion for the Sophist is now immediately recognized as the fundamental proportion of the Divided Line, which we formulate as a proposition:

\medskip

\noindent{\it Proposition.} \emph{The ratio of the second step in the Division of the Sophist, 
namely the ratio of the art of making real things $D_2$, to the art of making images $C_2$ 
is equal to 
the ratio of the fifth step in the Division of the Sophist,
namely the ratio of the mimetic art with knowledge $D_5$ to the mimetic art with opinion $C_5$.}

\medskip

 Plato takes extreme care to make sure that indeed the ratio of the Divided Line is applicable in the Logos of the Sophist:
 The ratio of real things B1 to images of real things B2 in the Divided Line is for all (sensible) things, both divine and human, while  the ratio of real things D2 to images of real things C2 in the Name+Logos of the \emph{Sophist} is only for human things.
Are the two ratios equal? We normally would equate the two ratios without thinking much; but Plato is more careful and presents a mathematically  structured argument in support of their equality. Here is the argument, in the \emph{Sophist} passage 265d-266d:

In the notation of the Divided Line, we have:

B1= the art of making real things;

B2 = the art of making images of real things.

In the notation of the division for the Sophist, we have:

C1= human productive arts;  

D1=divine productive arts; 

D2=the art of making human real things; 

C2= the art of making human real things.

Plato regards all these entities as line segments, and interprets D2 and C2 as rectangles:

D2= the art of making human real things as the rectangle B1$\cdot$C1, and  

C2= the art of making human real things as the rectangle B2$\cdot$C1.

He
then applies the philosophical analogue of Proposition VI.1 of the \emph{Elements} to obtain

B1/B2= B1$\cdot$C1/B2$\cdot$C1=D2/C2.

A similar (implied) argument yields 

A/B= D5/C5.

Now we consider the next two Logoi of the Name and Logos of the \emph{Sophist}, which are analogous to the geometrical ratios of periodic anthyphairesis.

The definition of the Sophist presents an additional feature that ties it even closer to the mathematical model: there are two Logoi after the Logos Criterion, and if the mathematical anthyphairetic model is indeed followed, we would expect to have two further equalities of Logoi, namely, the ratio  of the third division step should be equal to the ratio of the sixth division-step, and the logos of the fourth division-step should be equal to the ratio of the seventh and final division-step. 

Now we consider the equality of the third and sixth ratios of the Division of the Sophist:
 
In the third division-step, the Genus ``humanly produced images" is divided into two species, 
 the likeness-making/faithful to the original  (``eikastikes") and the fantastic/distorting the original (``phantastikes"). (236c3-7, 266d8-9);  
In the sixth division-step, described in 267e7-268a8, the Genus opinionated-imitator is divided into the two species, ``the simple" (``haploun") and ``the dissembling" (``eironikon") opinionated-imitator, 
Thus the simple imitators (step D6) do not distort (step D3) their opinion, but rather express a likeness of their opinion, while the dissemblers (C6) distort (C3) and disguise their opinion behind a false appearance. Therefore, we have D3/C3=D6/C6.

We go on with the equality of the fourth and seventh ratios:

In the fourth division-step, described in 267a1-b3,  in mimetic art the instrument of imitation is the imitator himself (C4), while in the nameless opposite art the instrument of imitation is other than the imitator (D4).

In the seventh division-step, described in 268a9-c4, he who listens to a dissembler is deceived and contradicted, 
if that dissembler is a demagogue, not by himself but by another instrument of deceit (namely the demagogue himself) (D7), 
while if the dissembler is a sophist, the listener is forced by the sophist to himself become the instrument of deceit (C7).

We conclude that D4/C4=D7/C7. 

The proportions D3/C3=D6/C6 and D4/C4=D7/C7 provide powerful additional evidence in favor of our interpretation of Plato's Name+Logos as a philosophical imitation of the geometric periodic anthyphairesis. 

Now we pass to the ``Name and Logos" of the \emph{Sophist}. 

\begin{figure}[!h] 
\centering
\includegraphics[width=1\linewidth]{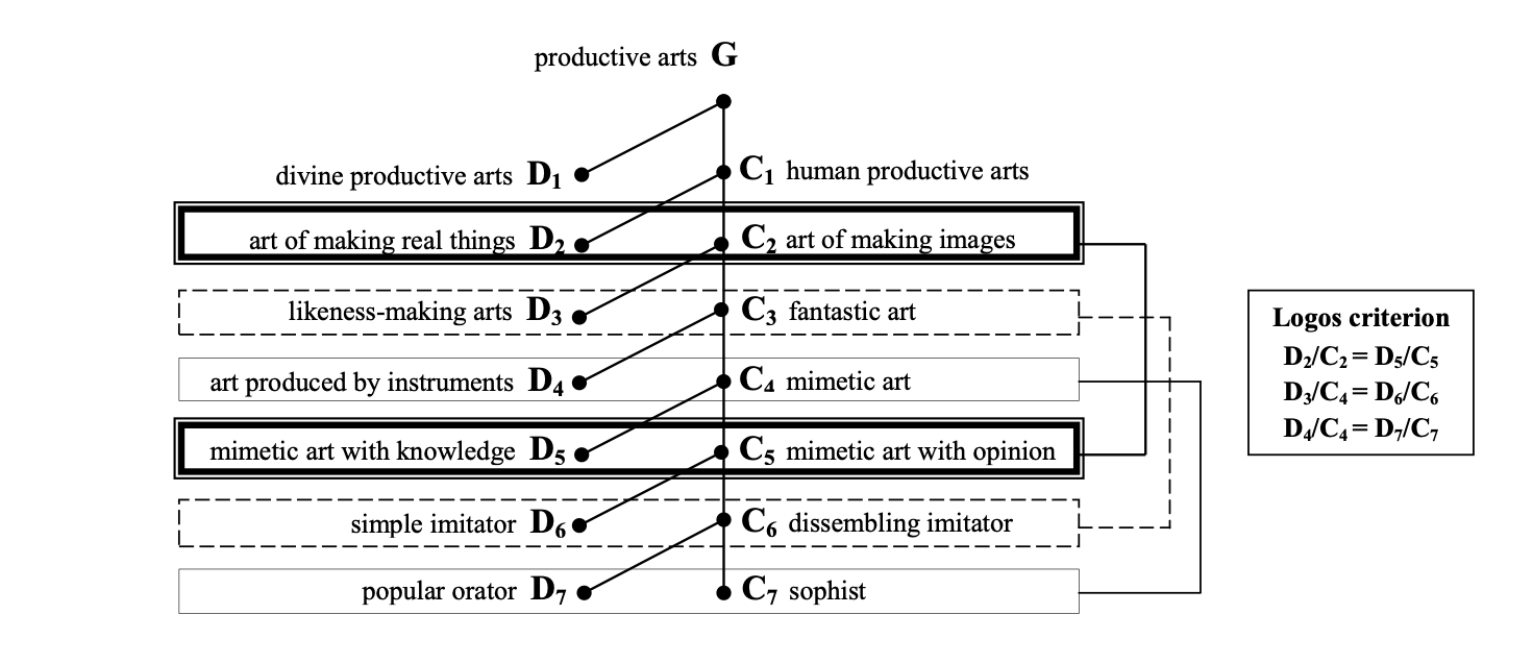}
\caption{\small {The definition of the Sophist by Name and Logos 
as the philosophic analogue of periodic anthyphairesis 264-268, 258b5-7}}
\label{fig:Sophist4}
\end{figure}

The complete Division and Collection of the Sophist (\emph{Sophist} 264b-268d) can be summarised in the scheme reproduced in Figure  \ref{fig:Sophist4}.

We recognize the proportion of the Divided Line D2/C2=D5/C5. Thus, we have full confirmation that Plato's Logos is the philosophical imitation of the geometric Logos Criterion (cf. \S 5).  Only this time, the Logos is not of the simple-minded, naive type that appeared in the definition of the Angler (from above downwards/from below upwards), but of a more sophisticated and philosophical kind, playing a central role in the dialectics of the Republic.

\section{The One of the second hypothesis  in the \emph{Parmenides} forms with the Being the philosophic analogue of the dyad in periodic anthyphairesis}\label{s:Parmenides}

Plato's \emph{Parmenides} is  considered to be a most difficult dialogue, with a plethora of apparent contradictions. We shall see in this section that the mathematics of periodic anthyphairesis provide the basis for an interpretation of the essential part of the dialogue, an interpretation which is in full agreement with the rest of Plato's dialogues.

The \emph{Parmenides} contains three main parts: the \emph{Introduction} (126a1-137c3), \emph{The One of the First hypothesis} (137c4-142a8)  and \emph{The One of the second hypothesis} (142b1-155e3).  
In this section we shall be mainly interested in the third part. Our goal  is to show  that the One of the second hypothesis is an intelligible Being, that its Knowledge  is described as Name plus Logos, and that it forms with its part Being a dyad which is the philosophical analogue of periodic anthyphairesis. 

Plato, in the second part of the \emph{Meno}, states the hypothesis: ``if the One is" (ἓν εἰ ἔστιν), 
meaning ``if the One is [an intelligible Being]" (142b1-5).  After an attempt he had made in the first part, with  a partless One, akin to a geometric point, to show that it is an intelligible being, which however failed, he introduces in the second part a second, highly non-partless One, radically different from the One of the first hypothesis, and he succeeds this time to show that it is an  Intelligible Being, as clearly stated in the conclusion of the second hypothesis (155d-e).

The One of the second hypothesis is  described as an entity that possesses Name and Logos (Καὶ ὄνομα δὴ καὶ λόγος ἔστιν αὐτῷ, καὶ ὀνομάζεται καὶ λέγεται, 155d8-e1).
In accordance with our previous interpretation of Name plus Logos, we expect that this One of the second hypothesis will have the structure of a dyad satisfying the philosophical analogue of periodic anthyphairesis. Indeed, the  first part of the second hypothesis (142b-143a5) is devoted to showing that there is an initial dyad, called (One, Being) satisfying the philosophical analogue  of periodic anthyphasirsesis.

Thus, let us see  how the 
dyad (One, Being) is formed, and the crucial role of the associated division. 
Our interpretation of the \emph{Parmenides} 142e4-5 will lead to the first two anthyphairetic relations of this dyad.

 The second hypothesis is interpreted as producing a dyad of parts, the One and the Being, in such a way that 
``the One participates in the Being" (οὐσίας μετέχει τὸ ἕν, 142b5-7, c5-7).
 The non-symmetrical expression ``the Being of the One" in 142b7-142b8 
(there is no corresponding expression like ``the One of the Being") suggests that 
the Being is a part of the One.
After some discussion, Plato declares that
 ``each of these two parts [of the One] possesses (ischei) the One and the Being" (142e3-4). We conclude from this that the One possesses two parts, the Being and another One, and that the Being also possesses a One and a Being. 
 These new parts will be generated from the two already existing ones, according to the following crucial rule:

καὶ γίγνεται τὸ ἐλάχιστον ἐκ δυοῖν αὖ μορίοιν τὸ μόριον (142e4-5).

This sentence is difficult to render precisely. Our rendering is that it is a rule of generating a new part (``to morion gignetai") from two already existing ones (``ek duoin moroin"), with the most crucial addition that the generated part will be the smallest (``to elachiston").
 
Our starting point is the dyad (One, Being). As noted above, the One contains as part the Being, and the two initial parts One and Being generate, according to the crucial rule, another part, say $x$, smaller than the parts that generate it, hence smaller than Being. But the One must possess two parts, a Being and a One, as noted, hence the part $x$ may be called $\mathrm{One}_1$. Thus, we obtain the equation

One=Being+$\mathrm{One}_1$, with Being$>$ $\mathrm{One}_1$,

\noindent clearly an anthyphairetic one.

Now the generation of the part $\mathrm{One}_1$ explains an earlier claim, in 
142d9-e3, that the Being contains a One part; indeed, the Being contains the part $\mathrm{One}_1$.
This sets forth the mechanism by the crucial rule, in 142e4-5, of the generation of the third smaller part from the two already existing ones, which implies that there is a Being part, say $\mathrm{Being}_1$, such that

Being=$\mathrm{One}_1$+$\mathrm{Being}_1$, with $\mathrm{One}_1>\mathrm{Being}_1$,

\noindent again an anthyphairetic relation. 

We have thus obtained the first two anthyphairetic relations of the anthyphairesis of the initial dyad (One, Being).
 
This anthyphairetic interpretation of the crucial \emph{Parmenides} 142d9-e5 passage (``part One")  is fully confirmed in a passage from Proclus' \emph{Comments on the Parmenides} 142d9-e5 which are contained in the \emph{Platonic Theology} 3, 89, 21-23.
Proclus writes: 
 
\begin{quote}\small

On the one hand, the One, by participating in (metechon) the Being, 
\\
is being divided (diaireitai) again, 
\\
so that the One and the Being generate (‘apogennan') 
\\
a secondary One (henada deuteron), 
\\
equal (suntattomenen) to a part (moirai) of the Being.

\end{quote}

Thus, Proclus' comment confirms our ongoing anthyphairetic interpretation, namely,

(i) Two parts One and Being ``generate" (``apogennan") a new part $\mathrm{One}_1$, as follows: the One contains the Being as part, and the two as a dyad generate a new partial $\mathrm{One}_1$ in such a way  that $\mathrm{One}_1$ is equal to a part of Being. 
In particular, Proclus' description fully agrees with our interpretation of the words ``gignetai ek" and ``elachiston".

(ii) The sentence ``the One participating (metechon) in Being is divided into One and Being" means that the One has parts in the sense that the One consists of parts Being and $\mathrm{One}_1$, with $\mathrm{One}_1$ equal to a part of Being.

Finally, we obtain the genuine meaning for the participation of the One into the Being: 
``The One participates in the Being" (in 142b5-7, c5-7)
means precisely that
``The One contains the Being as part, 
the (passive) One is divided by the (active) Being, and 
there is generated as anthyphairetic remainder $\mathrm{One}_1$ smaller than the Being".

We now turn to Proclus' revealing comments for the Division of the Being. The following is 
  in \emph{Platonic Theology} 3, 89, 24-26 and it concerns the part Being: 

\begin{quote}\small

On the other hand, the Being, by participating (metechon) in the One, 
\\
is again being divided (diakrinetai) in Being and One; 
\\
because it generates (apogennai) 
\\
a Being more partial (merikoteran) dependent on (exertemenon) 
\\
a more partial (merikoteras) One.

\end{quote}

It is notable that Proclus deals first with the division of the One, and
then, only as a by-product of this division, comes the generation of the $\mathrm{One}_1$,
equal to a part of Being, thus allowing the formation of the dyad $\mathrm{One}_1< \mathrm{Being}$, now in position to deal with the division of the Being. Thus, the
division of the Being does not take place independently of the division of
the One, but on the contrary depends on it.

One point needs to be clarified. This is the sentence ``the more partial one"
(‘henados merikoteras'), the same as the second one (``henada
deuteran"), which we have denoted by $\mathrm{One}_1$. Proclus probably considers this
obvious, since every new part must be generated. This is done in a linear
way, and, at the stage of the second relation, the only existing ones are
the original One and ``the second one", $\mathrm{One}_1$, which is indeed a more partial
one, since it is equal to a part of the Being, which is a part of the One.
Note also that Proclus later in the \emph{Platonic Theology} refers to the
division of the One in the more partial ones (``kermatizousa de to hen eis
tas merikoteras henadas", 4,80,4-5). Thus, the ``second one", $\mathrm{One}_1$, generated by
division, cannot but be one of the ``more partial ones". Thus, Being
participates in $\mathrm{One}_1$ and as a result Being is divided by $\mathrm{One}_1$  and
generates a more partial Being, $\mathrm{Being}_1$.

Thus, our interpretation is in agreement with Proclus' comments interpreting ``gignetai ek", ``to elachiston morion", ``ara" and  ``metechei" of the \emph{Parmenides} passage. Both the \emph{Parmenides} passage and Proclus' comments on it confirm that the \emph{Parmenides} passage 142d9-e5 describes the first two steps of the anthyphairesis of One to Being.

After describing the first double step of the division of the dyad (One, Being), the passage 142e5-7 continues with the description of the general step, rendered as follows
\begin{quote}\small

and thus always (``aei") in the same manner (as before) (kata ton auton logon) whatever becomes (``genetai") a part, 
always possesses (`ischei') a One and a Being 
because (``gar", 142e7) 
a (generated) One always (``aei")  possesses (``ischei") as part a Being and 
a (generated) Being always (``aei") possesses (``ischei") as part a One (142e5-7) 

\end{quote}

The expression ``and thus always in the same manner" (kata ton auton logon houtos aei, 142e5) 
refers to the immediately preceding statement, 
``the two parts generate a smallest part" (142e4-5) 
and turns this statement into a general rule of generation (by division-subtraction) of a third least part from a dyad of already generated unequal parts. 

Note also that in 142e5-7, the generalized statement 
``every part possesses a One and a Being"
holds because (``gar", 142e7) the generalized statement 
``every One possesses a  Being"
  and 
``every Being posseses a One"
holds, precisely in the same way as the initial statement 
``the One possesses a Being and the Being possesses a One"
holds,
``hence (`ara', 142e3) both the One and the Being possess a One and a Being". 

Thus the term ``gar" establishes precisely the same causal relation and has exactly the same force for the general case (142e5-7) than the term ``ara" (142e3) for the initial case (142d9-e5).

We are then led, reasoning by mathematical recursion, to the precise statement that for every natural number k,

\begin{quote}\small
the $\mathrm{One}_k$ possesses (``ischei") the $\mathrm{Being}_k$ and the $\mathrm{One}_{k+1}$, with the $\mathrm{One}_{k+1}$ is equal to a part of the $\mathrm{Being}_k$ , and 
\\
the $\mathrm{Being}_k$  possesses (``ischei") the $\mathrm{One}_{k+1}$ and the $\mathrm{Being}_{k+1}$, with the $\mathrm{Being}_{k+1}$ equal to a part of the $\mathrm{One}_{k+1}$.

\end{quote}

 The result of this inductive argument is that every part is \emph{ad infinitum} divided into two parts (142e7-143a2), hence the dyad (One, Being) is in fact infinite in multitude.\footnote{\label{n:multitude} An explanation of the word ``multitude" is necessary here, and it will be used especially in \S \ref{s:dialectic} about dialectic numbers and eristic numbers. In particular, we need stress the difference between ``multitude" and ``number",  in the sense that Plato uses these words.  The word ``multitude" is the translation of the Greek word plethos (πλήθος). 
For instance, Euclid defines in the \emph{Elements}  (Definition VII.2) a number as ``multitude of units" (plethos monadon, πληθος μονάδων), meaning a finite multitude of units.
Roughly speaking multitude can be thought of as cardinality. But
Plato distinguishes sharply between multitude and number. 
 In the \emph{Parmenides}, he speaks of the infinite multitude of the anthyphairetic remainders of the One to the Being.
In case we consider eristic numbers, there is no essential difference between multitude and eristic number. 
But  as we explain in the present section, the true number for Plato is the dialectic number, which sharply differentiates from multitude.}

  (143a2-3),
the philosophical anthyphairesis is infinite. Thus, the anthyphairetic division of the ``indefinite dyad"  (One, Being) yields

One=Being +$\mathrm{One}_1$, Being$>\mathrm{One}_1$,

Being=$\mathrm{One}_1$+$\mathrm{Being}_1$, $\mathrm{One}_1>\mathrm{One}_1$

\ldots

$\mathrm{One}_n$=$\mathrm{Being}_n$ +Onen+1, $\mathrm{Being}_n>\mathrm{One}_{n+1}$,

$\mathrm{Being}_n$=$\mathrm{One}_{n+1}+\mathrm{Being}_{n+1}$, $\mathrm{One}_{n+1}>\mathrm{Being}_{n+1}$,

\ldots \\
and there is then an \emph{infinite multitude} (``apeiron to plethos", 143a2) of remainders-parts of the anthyphairetic division

One$>$Being$>\mathrm{One}_1>\mathrm{being}_1>\ldots >\mathrm{One}_n >\mathrm{Being}_n >\ldots$.

 Next we expect that Plato will establish that the One satisfies not only Name (anthyphairetic division), but also Logos (the philosophical analogue of periodic anthyphairesis) as well. 
 This is done but in a quite cryptic and roundabout way, and in fact not in terms of Logos, but in terms of contact (Plato uses the word ``hapsis").

  The word contact (``hapsis") in the \emph{Parmenides} (the verb ``haptoito" in 138a3-7, ``hapsis" in 148d-149d) is the philosophical analogue of the ``ratio between two successive parts" in the sequence of remainders of the anthyphairesis of One to Being. 
Indeed, according to the passage 148d5-149d7:

-- a contact (``hapsis") is a relation between two consecutive (``ephexes") terms in the sequence of the anthyphairetic remainders of the One to Being (148e4-7, 149a4-6), ordered according to generation, and 

-- the contacts satisfy the general  
``plus one rule": 
\begin{quote}\small
number of terms =hapseis+1 (149a7-c3).
\end{quote}

It turns out that this  ``plus one rule"  in 148d5-149d7 originates in 
the ``plus one rule" for musical intervals:

\begin{quote}\small
number of consecutive terms=musical intervals +1.
\end{quote}

This will lead us to the periodicity of the anthyphairesis associated with the dyad (One, being).

So what exactly is an ``hapsis"? 

  It is crucial for an understanding of the One of the second hypothesis that the ``plus one rule" appears before Plato solely in one place: Pythagorean music. Once we realize that the only use for Plato of the concept of ``contact" is the enunciation of the ``plus one rule" in 149b1-c3, it becomes possible to reveal the nature of ``contacts-hapseis", because the ``plus one rule" is used elsewhere in Plato, in several passages on musical theory. Here is a list of some such occurrences:
  
(i) \emph{Republic} 546b5-6;  (ii) Proclus, \emph{Commentary on Plato's Republic} 2,36,21-25; 2,37,12-14;
(iii) \emph{Republic} 616d6-617d5; (iv) Proclus, \emph{Commentary on Plato's Timaeus} 2,237,3-15;
(v) Timaeus 34b4-35b6; (vi) Proclus, \emph{Commentary on Plato's Timaeus} 2,187,13-15; 2,188,6-9;\footnote{Here, Proclus has worked out Plato's system in exhaustive detail, and he has found that 
there are 34 terms (horoi) and 33 musical intervals-ratios} (vii) Timaeus Locrus, \emph{On Nature} (Peri phusios) 209,6-7;(viii) \emph{Timaeus} 36d2-4;(ix) Iamblichus, \emph{Theologoumena Arithmeticae} 48,10-14; 50,5-6;(x) Euclid(?), \emph{Sectio Canonis}, Proposition 9; (xi) Nicomachus, \emph{Harmonicum enchiridion} 12,1,32-40; (xii) Aristides Quintilianus, \emph{de Musica} 3,1, 32-82.

Since  a musical interval between two terms is just the ratio between the terms involved, we conclude that a contact/hapsis between two consecutive parts in the sequence of the parts of the (One, Being) can only be the ratio of these two parts.  
Thus, a contact (hapsis) in the One of the second hypothesis of the \emph{Parmenides} is 
 a ratio of two consecutive terms in the sequence of the anthyphairetic remainders 
of the dyad (One, Being). This provides us with an essential clarification of the word ``hapsis". 

We conclude that the rule enunciated in detail in 149b1-c3 between the number of consecutive terms and the connective hapseis is indeed the rule

\begin{quote}
number of consecutive terms = ratios of consecutive terms +1.
\end{quote}

It is remarkable that Plato, certainly on purpose, refers 
only to the number of terms, and never to the number of contacts.
He phrases the ``plus one rule"
not as ``number of terms = number of contacts +1", but  always as ``number of terms = contacts+1". 
This indicates that the units for the numbers are never the contacts, but solely the parts of the One.

Now we come to the most crucial point of the whole interpretation. 
In the passage 144c2-d4, Plato introduces a consideration of the property that ``the One is in the Being": since the Being is the Other of the One, the proposition stating that ``the One of the second hypothesis is in the Other" (145d6-e6) confirms that the One of the second hypothesis possesses this property. But why, and according to what criterion? The essential criterion for a One to be in the Other is well-hidden, in an earlier passage, where it was declared, on the basis of this criterion that ``the One of the first hypothesis is not in the Other". The criterion there is given, most surprisingly, in terms of contacts and their circularity (138a3-7). We have already clarified the meaning of a contact (hapsis) as a ratio of consecutive terms in the infinite anthyphairetic sequence of
consecutive terms in the infinite anthyphairetic sequence 
$$\mathrm{One}>\mathrm{Being}>\mathrm{One}_1>\mathrm{Being}_1>…>\mathrm{One}_k>\mathrm{Being}_k>\ldots$$
generated by the dyad (One, Being) of the second hypothesis; it follows that a circularity of contacts is exactly the Logos Criterion for anthyphairetic periodicity, expressed in a different language.

We now have all the pieces necessary for a satisfactory interpretation of the crucial passage 
138a3-7.
We shall show that the circularity (138a5-7) of contacts/hapseis (148d-149d),
interpreted as ratios of the terms of the sequence of successive remainders, 
leads to the Logos Criterion (One/Being=$\mathrm{One}_k/\mathrm{Being}_k$  for some $k$) for the anthyphairetic periodicity of the dyad (One, Being) of the Second Hypothesis.
 
In 138a3-7, we read  that the statement 
 ``The One is in the Other" (᾿Εν ἄλλῳ μὲν ὂν) 
is equivalent to the statement that
there are contacts/hapseis (haptoito) of many parts of the Other (pollachou) by many parts  (pollois) of the One (autou); furthermore
the contacts/hapseis should be in such a way that they form a circle (πολλαχῇ κύκλῳ ἅπτεσθαι).
This is how the circularity of contacts appears.

The proof of the anthyphairetic periodicity of the dyad (One, Being), outlined below,
results from a combination of the passages 142b-144d, 138a3-7, and 148d-149d:

(1) the One is [present] in the Being 
(144d2-4, 145b6-e6);

(2)  the One is in the parts of the Being
by means of its parts (144d2-4);

(3)  the One is in the other (138a3);

(4)  the Other of the One coincides with
the Being (143b); 

(5) hence, the statement ``the One is in the Being" in 144c-d
and the statement ``the One is in the Other" in 138a coincide.

(6)  Conclusion:
The One has contacts (138a, 148d-149d),
ratios between successive parts (148d-149d),
by means of
the many parts of the One (pollois 138a5, memerismenon, 144d2)
to the many places of the Being (pollachou, 138a5, 144d1),
forming a cycle (138a4, a6).

(7)  The contacts in the sequence $$\mathrm{One}>\mathrm{Being}>\mathrm{One}_1>\mathrm{Being}_1>\ldots > \mathrm{One}_k>\mathrm{Being}_k>\ldots$$
form a cycle, resulting in the 
anthyphairetic periodicity of the dyad (One, Being) of the second hypothesis of the \emph{Parmenides}.

Since the One of the second hypothesis explicitly satisfies the property that 
the One is in the other (\emph{Parmenides} 145b6-e6), we finally conclude, as expected, that
the dyad (One, Being) satisfies the Logos Criterion for the periodicity of the anthyphairesis of the One to the Being. 

Furthermore, the anthyphairesis of the dyad (One, Being) in the second hypothesis of the \emph{Parmenides} is not abbreviated, as in the Name plus Logos of the \emph{Sophist} (those of the Angler, and of the Sophist, as Intelligible Beings), but is complete. This further strengthens our anthyphairetic interpretation of Name plus Logos.

From \S \ref{s:Angler} to \ref{s:Parmenides}, we conclude that Plato's method of Division  and  Collection, equivalently Name and Logos, equivalently True Opinion and Logos, is 
not a Linnaeus classification scheme, 
as usually interpreted by modern Platonists,
but a philosophical imitation of periodic anthyphairesis.

A further strong confirmation of the periodic anthyphairesis
interpretation of the divisions in the \emph{Sophist}, in \S 9 and 7, is
provided by the circular description of the method:

\begin{quote}\small
σχεδὸν γὰρ αὐτὸν περιειλήφαμεν ἐν ἀμφιβληστρικῷ τινι τῶν ἐν τοῖς λόγοις 
περὶ τὰ τοιαῦτα ὀργάνων, ὥστε οὐκέτ' ἐκφεύξεται τόδε γε. 235b1-3
 
for we have almost got him into a kind of encircling net of the devices we
employ in arguments about such subjects, so that he will not now escape
the next thing.
 
ΞΕ.– Πάλιν τοίνυν ἐπιχειρῶμεν, σχίζοντες διχῇ τὸ προτεθὲν γένος,
πορεύεσθαι κατὰ τοὐπὶ δεξιὰ ἀεὶ μέρος τοῦ τμηθέντος, ἐχόμενοι τῆς τοῦ
σοφιστοῦ κοινωνίας, ἕως ἂν αὐτοῦ τὰ κοινὰ πάντα περιελόντες, τὴν οἰκείαν
λιπόντες φύσιν ἐπιδείξωμεν μάλιστα μὲν ἡμῖν αὐτοῖς, (264d10-265a1)

Stranger.– Then let us try again; let us divide in two the class we have
taken up for discussion, and proceed always by way of the right-hand part
of the thing divided, clinging close to the company to which the sophist
belongs, until, having circularly stripped him of all common properties
and left him only his own peculiar nature, we shall show him plainly
first.
\end{quote}

In fact, all the descriptions of the method of Division and Collection are
circular. Cf., in addition to the above two passages, \emph{Statesman} 273d5-e4, 285a4-b6,
286d8-287a6, \emph{Phaedrus} 266a3-6, the One of the second hypothesis in the
\emph{Parmenides} 144d4-e3, 148d-149d, 138a3-7 described in \S 7, and
\emph{Philebus} 15d4-8, they possess circularity, that is, periodicity.

\medskip

Equipped with this structure of an intelligible Being, we  provide definitive answers to fundamental questions, that were not be resolved by Platonists, concerning the following topics: the dialectic numbers, which are based on the anthyphairetic periodicity and the plus one rule, stating that the dialectic number of terms of a sequence is the (number of) ratios of successive terms plus one (stated in the \emph{Parmenides} 148d-149d); the description of the intelligible being as an Indivisible Line, a statement bordering on the contradictory; the  also seemingly contradictory \emph{Sophist}'s statement that ``the not-Being is a Being", based on the equalization of the two elements of the dyad defining an intelligible Being;  the more general self-similar Oneness of an intelligible Being, based on the equalization of all parts generated by the anthyphairetic division of an intelligible Being; and finally the Third Man Argument in the Introduction to the \emph{Parmenides}, appearing as a threat for Plato's theory, but essentially innocuous because of the self-similar Oneness.

 \section{Dialectic numbers, self-similar Oneness of the intelligible Being, and ``not-Being is a Being"} \label{s:dialectic} 
 
In this section, we elaborate on the role of dialectic numbers in an intelligible Being 
in terms of the Logos Criterion and the ``plus one rule" (148d-149d). We show that 
the indivisible lines coincide with the intelligible Being, that
the intelligible Being is both Infinite and Finite, and both One and Many (in a self-similar sense) (144d4-e3), and we explain
  the statement ``not-Being is Being" in the Sophist.
    
In the second hypothesis of Plato's \emph{Parmenides} there are two definitions of numbers: the eristic and the dialectic. Plato does not call them with these names, but he uses this distinction in several places (e.g. \emph{Philebus} 17a4) to indicate the distinction between the high philosophical (dialectic) and the low practice by the many (eristic).
After the proof, which we outlined in \S \ref{s:Parmenides}, that the dyad (One, Being) satisfies the philosophical analogue of infinite anthyphairesis, with the sequence of anthyphairetic remainders, infinite in multitude:
 
$$\mathrm{One}>\mathrm{Being}>\mathrm{One}_1>\mathrm{Being}_1>…>\mathrm{One}_n >\mathrm{Being}_n >\ldots,$$
 Plato suggests that the concept of number may be introduced in the intelligible Being One by taking as units exactly these remainders. Thus, e.g., the number three could be formed with the units $\mathrm{One}, \mathrm{Being}, \mathrm{One}_5$. It has not been realized by scholars of Plato that these numbers are eristic, exactly because they consist of unequal units, since the infinite sequence of anthyphairetic remainders is a strictly decreasing sequence. Anyone who reads the \emph{Philebus} 56c-d passage will realize that these cannot be numbers that Plato would accept. And indeed, he rejects them by first pointing out that they provide an infinite number of units (``pleista", 144c1-2), and, shortly after this, declaring that this infinity of units in an intelligible Being must be rejected (144d5-7).

 It turns out that it is possible and meaningful to introduce the good, dialectic definition of numbers only after the periodicity and not just after the infinity of the anthyphairesis of the dyad (One, Being) has been established (\S \ref{s:Parmenides} as well).  After this, the units will remain unequal, as they were before, but they become equalized by the following property, in the presence of periodic anthyphairesis: 
 
 \emph{The number of parts of the One is finite and equal to the number of parts of the Being (or in fact of the $\mathrm{One}_1$, of the $\mathrm{Being}_1$, or of any other anthyphairetic remainder).} 
 
This equality is the only criterion for the equalization of the part One with the part Being, or with the part $\mathrm{One}_1$, or with the part $\mathrm{Being}_1$, etc.

Let us see now how the ``plus one rule", described by Plato in 148d5-149d7 between terms and contacts/ratios and applied to the sequence of parts/anthyphairetic remainders of an intelligible Being generates the dialectic numbers.
 
We start with the definition of the dialectic numbers in an intelligible Being, in terms of the anthyphairetic periodicity and the ``plus one rule". In the passage \emph{Parmenides} 148d-149d  of 
  the One of the Second Hypothesis, the dialectic number, defined for a finite or infinite sequence of successive terms (of an interval) in the sequence of remainders
$$\mathrm{One}>\mathrm{Being}>\mathrm{One}_1>\mathrm{Being}_1>…>\mathrm{One}_n >\mathrm{Being}_n >\ldots,$$
is equal to 
(the multitude of all different ratios formed from successive terms in S) + 1.\footnote{Thus, the set (or the sequence ) of all the anthyphairetic remainders of the anthyphairesis of the One to the Being possesses 
an infinite multitude in the sense we explained this word in Footnote \ref{n:multitude}, but a finite dialectic number.
Plato avoids the use of the word ``number" even when counting contacts (hapseis, namely ratios of successive remainders), 
but he also abstains from  the use of the word ``multitude", and so he simply says: 
number of terms=contacts+1. }
 
For instance,  the dialectic number of the sequence of natural numbers 2,4,8,9 is equal to the number of different ratios $(2/4 =4/8, 8/9)+1$, that is, 3 and not 4. The dialectic number of the infinite sequence $(2,4,8, …,2n, 2(n+1),\ldots)$ is equal to 2.

Now in a periodic anthyphairesis, the sequence of successive remainders is periodic. Therefore the multitude of the different ratios formed from the whole infinite sequence of parts of the One is finite, and thus the whole dialectic number of an intelligible Being is finite, never exceeding the length of the anthyphairetic period plus one.

%
%

The combination of the statements obtained 
in 142b-143b (\S \ref{s:Parmenides}) and in 144d4-7 (present section) 
brings us to a rather difficult and seemingly contradictory  situation: 
the family of parts of the One
has been shown (in \S 10) to satisfy an infinity condition, namely, that 
``the parts of the One is infinite in multitude", while
``the parts of the One are finite in dialectic number".
The same distinction is explicit in the \emph{Sophist} 256e5-257a9.
In fact,  as we already stresses in Footnote \ref{n:multitude},
Plato draws a fundamental difference between ``(dialectic) number" and ``multitude" which incidentally is not understood by modern Platonists, who routinely translate ``multitude" by ``number".

%

Now the intelligible Being is both Infinite and Finite.  
It is Infinite since the anthyphairesis of the One to the Being is infinite; and
it is Finite, since the multitude of contacts is, by anthyphairetic periodicity, finite. Equivalently, because the whole number in an intelligible Being is finite.

The fact that the intelligible Being is a mixture of Infinite and Finite, which is stated in the \emph{Parmenides} for the One of the second Hypothesis, is also stated in the \emph{Philebus} 16c5-17a5 and 23c9-25e5. 
In the \emph{Parmenides} and in the \emph{Philebus} 16c5-17a5, the Infinite is the Infinite anthyphairesis, and the Finite is the finite length of the period, equivalently the finiteness of the whole number (as we just explained). In the \emph{Philebus} 23c9-25e5,
the Infinite is again an Infinite anthyphairesis/incommensurability, 
but the Finite is now the finite anthyphairesis/commensurability. 
The relation of the two finiteness  can be explained by Theaetetus' theorem: 
if $a,b$ are two lines incommensurable in line but commensurable in square 
(\emph{Philebus} 23c9-25e5 infinite and finite), 
then the anthyphairesis of $a$ to $b$ is eventually periodic 
(\emph{Parmenides}, \emph{Philebus} 16e5-17a5 infinite and finite) \cite{N1997, NFB2024}.

The next fact to know is that the indivisible line is identified with Plato's intelligible Being.
Indeed, Plato hints at the indivisible (line) at least twice, in the definition of the Noble Sophistry by Name and Logos in the \emph{Sophist} 226b1-231b8 (at 229d5-6) and in the \emph{Phaedrus} 277b5-c6 (at 277b7). The meaning of indivisibility is a mystery for modern Platonists, see \cite{Rashed}. According to Aristotle in  \emph{Metaphysics} 987b22,  ``Plato's kinds (ta eide) are numbers", thus these  numbers, called ``eidetic" (1086a5, 1086a8, 1088b34, 1090b25), are equalized (1081a25, 1083b24, 1091a25), hence they correspond exactly to what we have called dialectic numbers. The finiteness of these numbers implies that after the completion of the anthyphairetic period, no new units are generated, and in this sense an intelligible Being is indivisible. The complete argument leading to the identification of an indivisible line with the intelligible Being is given in \cite{N2024b}.

  The equality of the dialectic number of the parts of the One with the dialectic number of the parts of the Being is stated in the following passage of the \emph{Parmenides} 144d4-5:
 
 \begin{quote}\small
 
 — Καὶ μὴν τό γε μεριστὸν πολλὴ ἀνάγκη εἶναι τοσαῦτα ὅσαπερ μέρη. 
 
— ᾿Ανάγκη. 144d4-5

\end{quote}

 A usual but incorrect rendering of this passage, in the classical translation by H. N. Fowler, 1925, is ``And that which is divided into parts (‘to meriston') 
must certainly be as numerous as its parts (‘mere')",
as if ``the parts" (mere) were the parts of the ‘``meriston", 
turning the statement into a trivial tautology.
But according to the immediately previous passage 144d2-4 to 144d4-5,  ``the One is by necessity ‘memerismenon' (divided)",
because the One could not possibly be simultaneously wholly   present 
in all the ‘mere' (parts) of the Being part.
Hence the sentence 144d4-5 cannot be rendered as Fowler does. 

 Our rendering of the passage 144d4-5 is different:

It is clear that the ``meriston" of 144d4-5 is the ``memerismenon"  One of 144d2-4, while the ``mere" of 144d4-5 are the ``mere" of the Being of 144d2-4. Thus, the correct rendering of 14d4-5 is:
``The One is divided into as many (‘tosa') parts 
as (‘hosaper') the parts of the Being".

Plato employs this equality of  dialectic numbers as the criterion for the equalization of the One and the Being in the second hypothesis of the \emph{Parmenides}, an equalization that will lead us to
the self-similar Oneness of an intelligible Being, and an explanation of the fact,  that the not-Being is a Being, deduced from the \emph{Sophist}.

 The criterion by which the One and the Being are equalized is the fact, shown above, that the dialectic number of the parts of the One is finite and equal to the dialectic number of the parts of the Being. We quote Plato again:

\begin{quote}\small 
-- For  it [the Being]  is not divided, you see, into any more parts than one, 
but, as it seems, into the equal number as the One"144d7-e1

for the Being is not wanting to the One, nor the One to the Being, 
but being two they are equalized (exisousthon) throughout.

-- That is perfectly clear. (144e1-3)
\end{quote}

But it is clear that by the same criterion every part of the One of the second hypothesis is equalized to the One. Hence the paradigmatical intelligible Being, the one of the second hypothesis in the \emph{Parmenides}, is a One in a self-similar sense: the One is divided into an infinite multitude of parts, but all parts are equalized to each other.

Furthermore, the dyad of the paradigmatical intelligible Being in the \emph{Parmenides} is named (One, Being), but the general dyad of an intelligible Being in the \emph{Sophist} is named (Being, not-Being). With these names, the not-Being is equalized with the Being,   and is thus an intelligible Being as well. 
The statement ``not-Being is a Being", a cause of considerable confusion among Platonists, is thus naturally interpreted, with no difficulty. In \S 13 below, we shall outline its Zenonian origin.

As a conclusion to this section, and to make a link with modern mathematics, let us note the fact that it follows from our discussion that Plato, more than two millennia before the 19th-20th century set theorists, tried to give a definition and the construction of numbers out of some fundamental notion exists, and what could be better than the One?

\section{The Uniqueness of the intelligible Being and the failure of the Third Man Argument} \label{s:Third}

 In this section, we discuss the so-called \emph{Third Man Argument} in the \emph{Parmenides} 132a1-b2, its disturbing reconstruction by Gregory Vlastos (1954) \cite{Vlastos}, our analysis of the falsity of this reconstruction, and its resolution in terms of the self-similar Oneness of an intelligible Being.\footnote{Let us note that the expression \emph{Third Man Argument} is not due to Plato but  it appeared later, by Aristotle, presumably because it produces a third cause for any already found pair of participating and participated entities.}

Let us start by stating the passage on Third Man Argument in the \emph{Parmenides} 132a1-133a10.
\begin{quote}\small

 I fancy your reason for believing that each idea is one is something like this; 
when there are many things (polla) which seem to you to be great (megala),
you may think, as you look at them all, that there is one and the same idea (mia… idea he aute) over all of them (epi panta), and hence you think the great is one (hen to mega).

\ldots

But if with your mind's eye you regard the great itself (auto to mega) and the other great things (talla ta megala) in the same way, will not another great (hen ti au mega) appear beyond, by reason of which (hoi) all these (tauta panta) must appear to be great (megala)?

\ldots

That is, another idea of greatness (allo eidos megethous) will appear, in addition to (par') greatness itself (auto to megethos) and the objects which participate in it (ta metechonta autou) and another (heteron) again on all these (epi toutois \ldots pasin), 
by reason of which (hoi) all these (tauta panta) are great (megala); and each of your ideas will no longer be one, but will be infinite in multitude (apeira to plethos) 132a1-b2.
\end{quote}

%
%
%
 
 Thus, the Third Man Argument generates an infinite regress that puts into question the Platonic tenet that for every definite property there is a unique intelligible Being that is the cause of this property; according to the Third Man Argument there would not be one intelligible Being for Greatness, one for Beauty, one for Justice, and so on, but infinitely many. If the Third Man Argument (TMA)  holds, then it undermines the foundation of Plato's philosophy.
 Vlastos proposed a reconstruction of the Platonic TMA in the \emph{Parmenides} 132a1-b2 by introducing three hypotheses. One of them, which we refer to as the
Non-Identity (NI), says: ``If an intelligible Being participates in an intelligible Being, then it is not identical to it. In particular, an intelligible Being does not participate in itself."

Our objection is based on the fact that
the participation of an intelligible Being A in an intelligible Being B is modeled after ``the One participates in the Being", meaning, according to Proclus' revealing comments, and as we argued in \S \ref{s:Parmenides}, that A is divided by B and generates an anthyphairetic remainder C, in an intelligible anthyphairesis. This implies, according to our interpretation of intelligible self-similar Oneness, that the participated Being B is equalized to the participating Being A. This equalization implies that Vlastos' NI hypothesis fails, and with it the whole TMA. In fact if 
$$\mathrm{One}>\mathrm{Being}>\mathrm{One}_1>\mathrm{Being}_1>\ldots$$
is the infinite sequence of anthyphairetic remainders, and if $k$ is the least index such that, by anthyphairetic periodicity, we have $$\mathrm{One}/\mathrm{Being}=\mathrm{One}_k/\mathrm{Being}_k,$$ then 
the menacing infinite regress of the TMA reduces to an innocuous infinite sequence 
of mutually equalized Beings:
$$\Sigma_0 = \mathrm{One}>\Sigma_1 = \mathrm{One}_k>\Sigma_2 = \mathrm{One}_{2k+1}> \ldots >\Sigma_n = \mathrm{One}_{k+n-1}> \ldots$$
Thus the internal self-similar Oneness of an intelligible Being implies the (external) uniqueness, 
seemingly challenged by the Third Man argument.
The infinite regress 
$$\Sigma_0=\mathrm{One} >\Sigma_1= \mathrm{One}_k>$$
$$\Sigma_2=\mathrm{One}_{2k+1},\ldots>\ldots> \Sigma_n=\mathrm{One}_{k+n-1}>\ldots$$
 reduces to an infinite sequence of mutually equalized Beings.
Thus the internal self-similar Oneness of an intelligible Being implies the (external) uniqueness, seemingly challenged by the Third Man argument. 

In fact, Plato simply states TMA in the Introduction of the \emph{Parmenides} as a problem, and in effect implicitly deals with it in the \emph{Second hypothesis}, with the equalization of each of the parts of the One with the One itself, by means of anthyphairetic periodicity.
For more details, the reader may consult \cite{N2026}.

\medskip

The final part of our study aims to prove that, contrary to the presently dominant interpretation of Zeno's arguments and paradoxes as being devoid of mathematical content, the analysis of Zeno's presence in the \emph{Parmenides, Sophist} (via the Eleatic Stranger), and Zeno's verbatim Fragments preserved by Simplicius, show that Plato's intelligible Beings essentially coincide with Zeno's true Beings, and hence that Zeno's philosophical thought was already anthyphairetic, and hence heavily influenced by the Pythagorean's Mathematics. 
These findings run against Burkert's claim that ``ontology is prior to mathematics".

\section{The anthyphairetic, and eventually Pythagorean, nature of Zeno's arguments and paradoxes
}

 Plato's true or intelligible Beings satisfy the Compresence of almost contradictory properties, such as Finite and Infinite, One and Many, Motion and Rest.  But the notion of true Being was conceived before Plato, by Zeno of Elea. A crucial text of Zeno in this respect his Fragment B3, transmitted by Simplicius in his \emph{Commentary to Aristotle's Physics} 140,27-33 where Zeno writes verbatim:
\begin{quote}\small
If the Many are true Beings, 
then necessarily the Many are as many (tosauta) as (hosaper) they are, 
and neither more (oute pleiona) of them nor fewer (oute elattona). 
But if they are as many (tosauta) as they (hosa) are, they would be finite.
If the Many are true Beings, the Many are infinite beings.
For there are always other (hetera) in between (metaxu) beings, and again other (hetra) in between (metaxu) them. And in this way the beings will be infinite.
\end{quote}

In this section, we show how the Compresence of almost contradictory properties held by both Zeno's and Plato's true Beings, leads 
to a novel interpretation of the purpose of Zeno's paradoxes and arguments, namely to separate true Beings from sensibles, and
 to the identification of Zeno's true Being with Plato's intelligible Being.

  One central contradictory property is the fact that the One of the second hypothesis in the \emph{Parmenides}, the paradigmatical intelligible Being, is Infinite and Finite. We have interpreted this by the fact that it possesses infinite anthyphairesis, and at the same time it is
  Finite, in the sense that it possesses periodic anthyphairesis, or, equivalently, because its whole dialectic number is Finite, or, equivalently, because the ratios of successive anthyphairetic remainders of the (One, Being) is Finite (\emph{Parmenides} 145a2-3).  Another property is that the One of the second hypothesis is 
Many because the One is divided, by anthyphairetic division, into an infinite multitude of parts and it is One because all these parts are equalized in the sense that they possess equal finite number of parts (\emph{Parmenides} 145a2-3).

Another remarkable property of Plato's 
intelligible Being is that it is 
in Motion and at Rest (\emph{Parmenides} 145e7-a8).

%

Zeno, in the passage \emph{Parmenides} 128e6–130a2, sets the question to Socrates whether true Beings satisfy the Compresence of Opposites properties: 
similar and dissimilar, one and many, motion and rest, is closely connected with infinite and finite, cf. \emph{Parmenides} 158e1-159a6. This Compresence  may be accepted with utmost difficulty, because it is not only that these entities are separately, in different ways, P and not-P, but also that the intelligible Being with P is at the same time an entity with not-P and P. Opposite properties for sensibles do not participate to each other (e.g. sensible motion and sensible rest). But intelligible true Beings do participate, e.g. the same/tauton and the other/heteron.

Most Platonic scholars believe that Zeno's suggestion is eventually rejected, but this clearly runs against the sense of Plato's \emph{Parmenides} text, namely, the fact that Plato's intelligible Beings definitely satisfy these properties, with carefully constructed proofs in the second hypothesis of the \emph{Parmenides}, and also because of Proclus'  clarifying comments, which we review now.  

Indeed, Proclus, in his \emph{Commentary to the Parmenides} 757,15–758,20, carefully interprets the \emph{Parmenides} 128e6–130a2 passage as a gradual shift and conversion of Socrates, gently guided by Zeno, from rejection, to neutrality, and then to enthusiastic acceptance of the Compresence of Opposites in the true Beings:

--Stage 1 of Socrates' conversion: 
from rejection (apognosis) to ``monstrosity", with the description of this mixture as something monstrous (teras);

– Stage 2: neutral ``amazement",
to a state of neutrality (huponoian \ldots tou alethous), with the description of the mixture as something causing puzzlement (thaumasomai, axion thaumazein);

– Stage 3: acceptance, ``admiration''
and finally positive certainty (dia tes apodeixeos bebaiothen, he teleutaia psephos alethestate), with the description now as definite enthusiasm (agaimen, thaumastos, andreios, agastheien).

 We are led now to our novel interpretation of Zeno's arguments and paradoxes: Plato's purpose of establishing the Compresence of almost contradictory properties in his true Beings is meant to separate the true Beings from the Sensibles.

  To achieve this interpretation, we first need to show that Zeno, by the Many, means the Sensibles.
  
This is explicitly supported in the crucial passage \emph{Parmenides} 128e6–130a6, in which Socrates refers systematically to ``the many" (ta polla) as ``the sensibles": (a) in   129a2–3 (``and that you and I and all others (talla) which we call many (polla) participate (metalambanein) in these two (toutoin duoin ontoin) [the idea of likeness and the idea of unlikeness]; (b) again, in 129d3–4 (``the many" include ``stones (lithous), sticks (xula), and the like"); and, (c) finally,  in 130a1, as in the visible objects (en tois horomenois).

Furthermore, Simplicius, in his \emph{Commentary to Aristotle's Physics} 97,9–99,6, quotes a
passage of Eudemus, where ``the Sensibles" are expressly called ``the Many" by Zeno:
\begin{quote}\small
 [Zeno] was puzzled (eporei) as it seems because each sensible (ton men aistheton  hekaston) was said to be many (polla) both categorically (kategorikos) and by means of division (merismoi), 97,13–15).
 \end{quote}

Plato, in the \emph{Republic} 475c–480a, systematically calls the Sensibles that participate in the intelligible Being F ``the Many F". In Plato's \emph{Parmenides passage} 132a1-b2, describing the Third Man Argument, 
``poll"  (many) in 132a2 has the meaning of ``sensibles".

From this we conclude that in Zeno's terminology, ``the many" are ``the sensibles". 

This will be useful in arriving at the correct interpretation of Zeno's fundamental hypothesis ``ei polla estin" in the immediately preceding passage \emph{Parmenides} 127e–128a (Zeno's Basic Argument), and in Zeno's Fragments B1 and B3. 

Now we come to our interpretation of Zeno's Basic Argument (\emph{Parmenides} 127e–128a).
We obtain this interpretation 
in five steps:
 
[Step 1] We assume that the Many are (ta polla einai), namely that a Sensible is a true Being.

[Step 2] A true Being satisfies one of the Compresence of Opposites property.

[Step 3] Then the Sensible must satisfy one of the Compresence of Opposites property.

[Step 4] But from our experience, we know that this is impossible, a contradiction.

[Step 5] Hence it is not true that the Many are, namely a Sensible is different from a true Being.

This provides a novel interpretation of the form and the purpose of Zeno's paradoxes and arguments: 
We claim that the purpose of these paradoxes and arguments  is to help Zeno's  teacher, Parmenides, against the attacks that he suffered, but not by showing that motion and multiplicity are impossible, but rather by showing that the sensibles are different from the true beings. We take as example his third paradox of motion, the so-called arrow paradox:

``The arrow is at rest as it moves" (Aristotle's \emph{Physics} 239b30).

This is not an argument against motion, but an argument for the difference between
intelligible and sensible motion, and it should be read as follows:

[Step 1] Suppose that the sensible is a true Being.

[Step 2] A true Being is at rest as it moves.

[Step 3] Then the arrow, a sensible entity, must be at rest as it moves.

[Step 4] But from our experience, we know that this is impossible, a contradiction.

[Step 5] Hence motion of a true Being is different from sensible motion, and the sensibles are
different from true Beings.

The argument with which Step 2 is established in Plato's \emph{Parmenides} 145e7–146a8 and \emph{Sophist} 248d10–250c5 is essentially the same as Zeno's argument in Aristotle's \emph{Physics} 239b6. That a sensible cannot be at rest as it moves [Step 4] is shown in the \emph{Republic} 436c5–437a3.

Furthermore, Zeno's Fragment B3 compared with the One of the second hypothesis in Plato's \emph{Parmenides} reveals that Zeno's true Being has the structure of a periodic anthyphairesis. Our interpretation of Zeno's Basic Argument applied to Zeno's Fragment B3 which we quoted at the beginning of this section, becomes:

[Step 1] Suppose that the sensible is a true Being.

[Step 2] A true Being is infinite (in multitude) and finite (in number).

[Step 3] Then a sensible entity must be infinite and finite.

[Step 4] But from our experience, we know that this is impossible, a contradiction.

[Step 5] Hence a true Being is separate and different from the sensibles.

  Step 2 should be put in parallel with Plato's \emph{Parmenides} 145a2–5 and \emph{Sophist} 256e5–257a12. 

Let us highlight now the close connection between Zeno's Fragment B3 and the periodic anthyphairesis of the One of the second hypothesis in the Parmenides.

\noindent {\it Infinity.} First of all we clarify the meaning of ``metaxu". This word, in Fragment B3, is almost always translated incorrectly as ``between" the two. But it should be rendered as ``offspring of", ``generated by" (ekgonon) the two, as explicitly noted by Simplicius by the word ``metaxu", not in the sense of separating both, but in the sense of consisting of both: μεταξὺ δὲ οὐ τὰ κεχωρισμένα ἀμφοῖν, ἀλλὰ τὰ ἐξ ἀμφοῖν συνεστῶτα (186,11), cf. Plato's \emph{Timaeus} 50d2–4 

Infinity of Zeno's true Being is established in Fragment B3 by the repeated \emph{ad infinitum} generation of parts by ``heteron". But in the One of the second hypothesis of the \emph{Parmenides}, the expression ``the One is heteron to the Being" is identified with the expression ``the One participates in the Being". And as we saw, participation was interpreted by Proclus to lead to the philosophical analogue of the infinite anthyphairesis of the dyad (One, Being). We conclude that the infinite generated in Zeno's B3 is closely associated with the anthyphairetic infinite generated by the One of the second hypothesis in the \emph{Parmenides}. 

\noindent {\it Finiteness and Equalization.} Finiteness of Zeno's true Being in Fragment B3 is established by the claim ``so many… as many" (``tosa…hosaper") and therefore finite. But in \emph{Parmenides}' One of the second hypothesis, the exact same expression ``tosa…hosaper" was interpreted to mean the periodicity of the anthyphairesis of the One and Being, and the subsequent equalization of the two parts One and Being. Thus, the Finite in Zeno's B3 is closely associated with the anthyphairetic periodicity by the One of the second hypothesis in the \emph{Parmenides}. 
With these remarks, Zeno's Fragment B3 has an anthyphairetic content and shows that Zeno's philosophical thinking is anthyphairetic.

We outline another argument showing the close similarity between Zeno's true Being and Plato's intelligible Beings

In the pseudo-Aristotelian treatise \emph{On Indivisible Lines} 968a18-b4, the author is supposed to show that there is an entity initially described as ``partless magnitude" (megethos ameres, 968a19, a22), and as ``indivisible line" (atomos gramme, 968b4). Furthermore, an indivisible line is identified with a ``partless magnitude", a clear variation of the ``magnitude-less magnitude" description of Zeno's true Being in Fragment B1. We conclude that Zeno's true Being coincides with the indivisible line.

We now outline the strong evidence that not only Plato's Compresence of almost contradictory properties, but also Plato's anthyphairetic statement ``not-Being is a being" that essentially produces them, is due to Zeno.

(1) Parmenides sternly warning against ``not-Being is" in his Fragment 7.1, mentioned twice in the \emph{Sophist} 237a8-9, 258c6-d4; 

(2) The Eleatic Stranger in the \emph{Sophist} 241d3 is concerned that the philosopher who introduced not-Being would be considered, by introducing the not-Being, as guilty of  ``parricide against his father [Parmenides]";

(3) Simplicius in his \emph{Commentary to Aristotle's Physics} 138, 22-25, that Zeno introduced not-Being in order to prove the One and Many property for true Beings;

(4) Parmenides shows his displeasure and opposition to Zeno's lecture on One and Many (Introduction of the \emph{Parmenides}, 127c5-d5); 

(5) In Plato's account of ``the battle of the Giants" (Sophist 245e6-247c8), an unnamed philosopher introduces a definition of true Being of the type ``not-Being is a Being" proving the One and Many Compresence property;

(6) Zeno confirmed, in the Introduction of the \emph{Parmenides} 128a-b, that Parmenides proved that the true Being is ``One" (128b3) [only One], while Zeno proved that it is ``not Many" (128b4) [not only Many]; and, 
 
(7) Zeno was clearly apologetic about composing his arguments [on Compresence]: he was very young, his manuscript was not intended for publication, but it was stolen (introduction of the \emph{Parmenides}, 128d6-e1). 

We can confidently conclude that Zeno was the person who introduced not-Being, namely the anthyphairetic dyad in periodic anthyphairesis as his true Being. 

Zeno's philosophical thinking can be anthyphairetic only if it were influenced by the early Pythagoreans, Pythagoras or more likely Hippasus. We can describe briefly the relation as follows: the early Pythagoreans, most likely Hippasus, proved the incommensurability of the diameter to the side of a square by computing the infinite anthyphairesis of the diameter to the side of a square \cite{NF2025}. This great mathematical achievement by the early Pythagoreans was exploited by Zeno to conceive of a true Being with the structure of periodic anthyphairesis, thus introducing the ``not-Being", a notion to which the philosopher Parmenides had strong objections, but later adopted by Plato, even though the ``philosophical children" of Parrmenides, Zeno/Eleatic Stranger in the Sophist, might be accused of ``parricide".

The anthyphairetic nature of Zeno's true Being is the absolute argument (a) against Burkert's, 1962/1972 thesis that ``ontology is prior to mathematics";  in fact Zeno was most heavily influenced by the early Pythagoreans, and (b) against Netz's \cite{Netz1999} and \cite{Netz2022} thesis that ``the mathematician Pythagoras perished finally at 1962 AD" and that ``most important: Pythagoras the mathematician was, I argue, indeed a myth". Both Burkert's and Netz's claims are thus finally rejected. Paraphrasing a famous misquote of Mark Twain, Pythagoras the mathematician, upon reading it, was heard saying that ``The reports of my death are greatly exaggerated" \ldots

Modern Platonists have never obtained a clear description of the structure of an intelligible Idea in terms of the mathematics of periodic anthyphairesis, and thus were not able to answer fundamental questions, nor to realize the close connection between Plato's intelligible beings with Zeno's true Beings.

To conclude, we come back to the last question in the list of 18 questions that we addressed in the Introduction, for which we venture an explanation: Plato's philosophy was the first and for many centuries the only one based on the mathematically inspired Finitization of the Infinite, and as such it was a hopefully consistent approximation of the contradictory, and therefore a system with great deductive power \cite{FN-Turaev}.  The two real architects of the system, Zeno and Plato,  consciously employed the paradoxical nature of the system in order to surprise their audience as much as possible \cite[\S 6.3.2]{N-Zeno}.

 \bigskip

\noindent{\bf Authors' addresses:}

\noindent Stelios Negrepontis, 
National and Kapodistrian University of Athens,
Mathematics Department, 
Panepistimioupolis, 15784, Athens, Greece.

 \noindent email:  snegrep@math.uoa.gr

\medskip

\noindent Athanase Papadopoulos,
Institut de Recherche Mathématique Avancée
(Université de Strasbourg et CNRS),
7 rue René Descartes,
67084 Strasbourg Cedex France.

\noindent email: papadop@math.unistra.fr

 \end{document}